\documentclass[10pt,leqno]{amsart}
\usepackage[pdftex]{graphics}
\usepackage{amssymb,amsfonts,amsmath,amsthm,latexsym,verbatim,epsfig}

\newcommand\beque{\begin{equation*}}
\newcommand\beq{\begin{equation}}
\newcommand\eeq{\end{equation}}
\newcommand\eeque{\end{equation*}}
\newcommand\eqn{\begin{eqnarray}}
\newcommand\beqna{\begin{eqnarray*}}
\newcommand\eeqna{\end{eqnarray*}}
\newcommand\feqn{\end{eqnarray}}

\newcommand{\nn}{\nonumber}

\newtheorem{lemma}{Lemma}
\newtheorem{propo}{Proposition}
\newtheorem{theorem}{Theorem}

\newtheorem{open}{Open Problem}
\newtheorem*{open2}{Open Problem}
\newtheorem*{theorem3}{Theorem}
\newtheorem*{lemma2}{Lemma}
\newtheorem*{propo2}{Proposition}

\numberwithin{equation}{section}

\begin{document}


\title[]{\Large{Error Estimates in Horocycle Averages Asymptotics: \\ Challenges from String Theory}}

\author{Matteo A. Cardella}

\address{Institute of Theoretical Physics
University of Amsterdam.
Science Park 904
Postbus 94485
1090 GL
Amsterdam,
The Netherlands.}

\email{matteo@phys.huji.ac.il}

\begin{abstract}
There is an intriguing connection between the dynamics of the horocycle flow in 
the modular surface  $SL_{2}(\pmb{Z}) \backslash SL_{2}(\pmb{R})$ and the Riemann hypothesis.
It appears in the error term for the asymptotic of the horocycle average of a modular function
of rapid decay.
We study whether similar results occur for a broader class of modular functions, 
including functions of polynomial growth, and of exponential growth at the cusp.
Hints on their  long horocycle  average are derived by translating the  horocycle
flow dynamical problem in string theory language. Results are then proved by
designing an unfolding trick involving a Theta series, related to the spectral Eisenstein
 series by Mellin integral transform.
We discuss how the string theory point of view leads to
an interesting open question, regarding the behavior of long horocycle
averages of a certain class of automorphic forms of exponential
growth at the cusp. 

\end{abstract}

\maketitle

\section{\Large{\textbf{Introduction}}}

\vspace{.1 cm}

In this paper we exploit  a novel angle for obtaining some insights on the long horocycle average asymptotic for
 certain classes of $SL_{2}(\pmb{Z})$-invariant  automorphic functions.
    We focus on  modular functions  of
 polynomial growth at the cusp, and on  a certain class of modular functions of (bounded) exponential growth.
 Automorphic  functions with such growing conditions  play a role  in string theory,  in the context of  perturbative (one-loop) closed string amplitudes.
Remarkably, their  horocycle averages
contain information on the   numbers of physical degrees of freedom
         of  closed strings particle-like excitations\footnote{A genus one closed string vacuum amplitude $\mathcal{A}$ is given by the integral of
a $SL_{2}(\pmb{Z})$ invariant function $f$ on the fundamental domain $\pmb{\mathcal{D}} \simeq SL_{2}(\pmb{Z})\backslash \pmb{\mathcal{H}}$,
$\mathcal{A} = \int_{\pmb{\mathcal{D}}} d\mu f$.
The effective numbers of closed string  states are encoded in the expansion  $ \int_{0}^{1}dx f(x,y) = \sum_{n = 0}^{\infty}(d_{n}^{B} - d_{n}^{F})e^{-\pi m^{2}_{n}y}$  of the automorphic function $f$ horocycle average, where $d_{n}^{B}$($d_{n}^{F}$) is the number of bosonic(fermionic) 
degrees of freedom  at    mass level $m_n$, $n = 1,2,\dots $.  Convergence of the long horocycle limit $y \rightarrow 0$ corresponds to
     a subtle pairing  among bosonic and fermionic closed string physical degrees of freedom. This cancelation was  called asymptotic supersymmetry in [KS]. Quite interestingly, in closed string theory  horocycle averages  asymptotics as (\ref{theta2}) when translated in closed  show an intriguing
         relation between  asymptotic supersymmetry and  the Riemann hypothesis [C1], [CC1], [CC3], [ACER].} [C1],[CC1],[CC3],[ACER].

\vspace{.35 cm}

The advantage  of  translating the dynamical problem in string theory language  is in the
possibility of using  consistency conditions from
string theory to gain insights  on the horocycle average asymptotic. For the  two classes of modular forms we  focus on,
the string theory perspective suggests a universal  behavior of their long horocycle average,
which  appears somehow surprising from the perspective of the theory of automorphic forms.
Our results are  then obtained by   an unfolding method
that  involves a  Theta series,  connected to the spectral Eisenstein series
by Mellin integral transform.  We
 illustrate advantages of the  Theta unfolding  for dealing with automorphic forms of not so mild growing conditions,
 over the classical Rankin-Selberg method.
In particular, we derive some  results previously obtained  by Zagier [Za2] via considerably
shorter proofs on the  analytic continuation of the Rankin-Selberg
integral transform for automorphic functions of polynomial growth.
 We then obtain asymptotics for long horocycle averages of modular functions of polynomial growth,
 including a relation between  the error estimate and  the Riemann hypothesis.
 For modular function of rapid decay the same kind of relation  was originally obtained  in [Za1].

When applied to modular functions playing a role  in string theory, our results
lead to fascinating  connections between enumerative properties
of closed string spectra and the Riemann hypothesis [C1],[CC1],[ACER],[CC2],[CC3].
These connections extend to multi-loops closed string amplitudes [CC1],[CC3]
and results for measure rigidity of unipotent flows   in homogenous spaces [Ra] are intertwined with properties of perturbative
closed string theory [CC2].

\vspace{.5 cm}

Let $\pmb{\mathcal{H}} = \{z = x + iy \in \pmb{\mathbb{C}}| y > 0 \}$ be the upper complex plane,
 horocycles in $\pmb{\mathcal{H}}$ are both circles tangent to the real axis in
rational points (cusps),  and  horizonal lines, (which can be thought as circles tangent to the $z = i\infty$ cusp).

 $\begin{pmatrix}a & b \\ c & d \end{pmatrix} \in SL_{2}(\pmb{R})$ acts on   $z \in \pmb{\mathcal{H}}$ through the  M\"obius transformation
$z \rightarrow \frac{az + b}{cz + d}$.
The following one-parameter  action of the upper triangular  unipotent subgroup $\pmb{U} \subset SL_{2}(\pmb{R})$

\vspace{.27 cm}

 \beq
 \pmb{g}_{u} : =
\Big\{
  \begin{pmatrix}
  1 & t \\
  0 & 1
  \end{pmatrix}, \,
    |t| \le |u|,     \,  \Big \}, \qquad u \in \mathbb{R} \nn
  \eeq

\vspace{.27 cm}

generates motions along  horizontal lines in $\pmb{\mathcal{H}}$. Long horocycles   in $\pmb{\mathcal{H}}$ do not exhibit  interesting dynamics in the half-plane $\pmb{\mathcal{H}}$, since
  the orbit $\pmb{g}_{u}(x + iy ) = \{ x + t + iy,  |t| \le |u| \}$ for $u \rightarrow \pm \infty$ just escapes to infinity.
However, $\pmb{g}_{u}(x + iy )$ has an interesting dynamics in the quotient space $\pmb{\Gamma} \backslash \pmb{\mathcal{H}}$, $\pmb{\Gamma} \simeq SL_{2}(\pmb{Z})$.  The horocycle $\pmb{g}_{u = 1}(x + iy)$ is a closed orbit
  in  $\pmb{\Gamma} \backslash  \pmb{\mathcal{H}}$
  with length   $1/y$, as measured by the hyperbolic metric $ds^2 = y^{-2}(dx^2 + dy^2 )$. Quite remarkably,  in the long length limit $y \rightarrow 0$, the horocycle
$\pmb{g}_{u = 1}(x + iy)$ tends to cover  uniformly the modular domain  $ \pmb{\Gamma} \backslash \pmb{\mathcal{H}}$ [He], [Fu],[DS],
(see figure 1. for plots of horocycles in the modular domain of increasing length obtained with Mathematica)

\begin{figure}[htbp]
\centering
\includegraphics[height = 12cm, width =  2cm ]{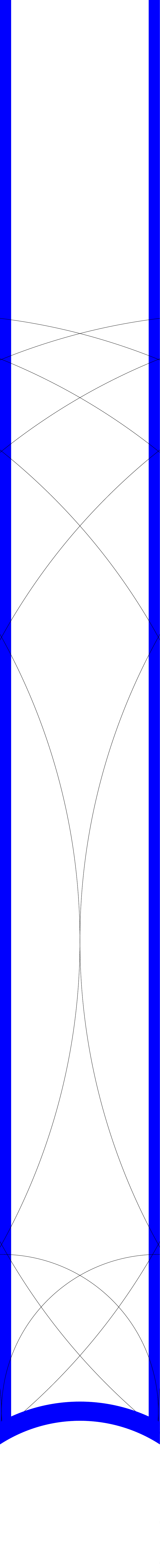}
\hspace{.8 cm}
\includegraphics[height = 12cm, width =  2cm ]{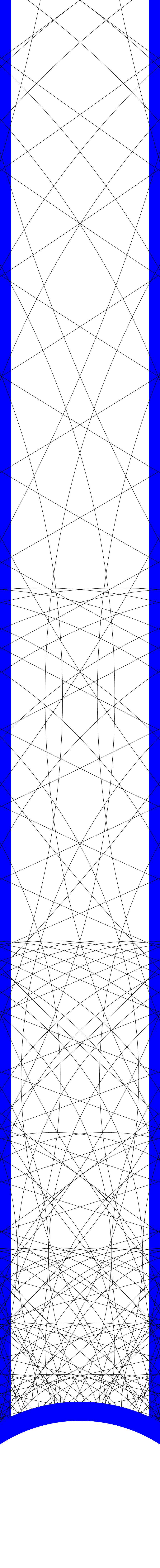}
\hspace{.8 cm}
\includegraphics[height = 12cm, width =  2cm ]{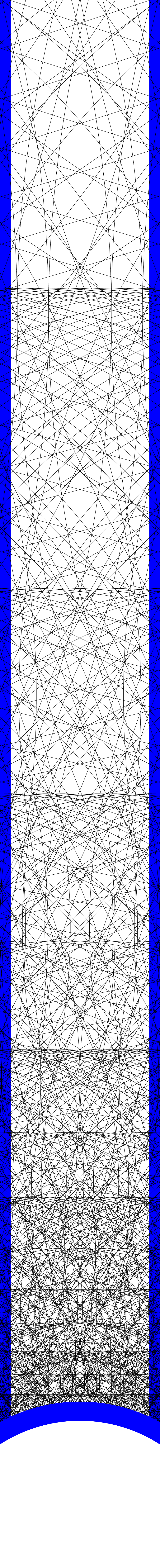}
\caption{ \textbf{Modular images  of  horocycles of increasing length.}
It is interesting to study the image  of $H_{\alpha} = \mathbb{R} + i\alpha$
in the standard  $SL_{2}(\pmb{Z})$ fundamental domain as $\alpha \rightarrow 0$.
Left:  modular image of the  line $y = \frac{1}{8}$.
Center: modular image of the line  $y = \frac{1}{100}$. Right: modular image
of the line  $y = \frac{1}{400}$.
In all cases the modular domain is truncated 
at $y > 10$.   The modular image of a line $y = \alpha$ tends to become dense for $\alpha \rightarrow 0$ [He].}
\label{Sflow}
\end{figure}

 Methods involving the theory of automorphic forms lead to  interesting results for  horocycle flow asymptotic.
Quite remarkable  is   the relation between error estimates for asymptotics involving
the average of an  automorphic forms  along long horocycles and the Riemann hypothesis.
By using the  Rankin-Selberg method, Zagier [Za1] has obtained the intriguing  result

\vspace{.27 cm}
\beq
  \int_{0}^{1}dx f(x,y) \sim \frac{3}{\pi}\int_{\pmb{\mathcal{D}}}d\mu f + O(y^{1 - \frac{\Theta}{2}}),  \qquad y \rightarrow 0 \label{theta2}
  \eeq

\vspace{.27 cm}

  when $f$ is a smooth   modular invariant function
   of rapid decay at the cusp  $y \rightarrow \infty$. Indeed, in order to have a sufficient condition for (\ref{theta2})
    to hold, one has to add some smoothness condition on $f$
    and a growing condition on  its Laplacian  $\Delta f$, [Ve],  (this is discussed in details in \S \ref{SPrapiddecay}, proposition \ref{p-rapid}).
   In eq. (\ref{theta2}) $\mu$ is the hyperbolic $\pmb{\mathcal{H}}$ measure, $d\mu = y^{-2}dx dy$, and
   the error estimate is governed by $\Theta = \rm{Sup}\{\Re(\rho) |  \pmb{\zeta}^{*}(\rho) = 0 \}$,  the superior of the real part of the non trivial zeros of the Riemann zeta function $\pmb{\zeta}(s)$, ($\pmb{\zeta}^{*}(s) = \pi^{-s/2}\pmb{\Gamma}(s/2)\pmb{\zeta}(s)$).

  \vspace{.27 cm}

   The error estimate for the convergence rate  in (\ref{theta2})  is remarkably linked to the   Riemann hypothesis (RH)\footnote{See also [Sa], [Ve] for a study of convergence rates for  horocycle flows and Eisenstein series for more general  quotients  $\pmb{\Gamma} \backslash \pmb{\mathcal{H}}$, where $\pmb{\Gamma} \subset SL_{2}(\pmb{Z})$ is a lattice.},
   indeed, RH is equivalent  to  the following condition

\vspace{.27 cm}

\beq
  \int_{0}^{1}dx f(x,y) \sim \frac{3}{\pi}\int_{\pmb{\mathcal{D}}}d\mu f + O(y^{3/4 - \epsilon}),  \qquad y \rightarrow 0 \label{Riemann}
  \eeq

\vspace{.27 cm}

   for every $f \in  C_{00}^{\infty}(\pmb{\Gamma} \backslash \mathcal{H})$.
Up to date,  the error term  is $o(y^{1/2})$ unconditionally  as a consequence  of the bound $\Theta \le 1$
   on the real part of  the Riemann zeta functions  zeros $\rho$'s.

\newpage

\centerline{\Large{\textbf{Notation and Terminology}}}

\vspace{.3 cm}

$\bullet$  $\pmb{\mathcal{H}} = \{ z = x + iy \in \mathbb{C}, \, \, y > 0    \}$, the upper complex plane.
\vspace{.2 cm}

$\bullet$ $\pmb{\Gamma}  \simeq SL_{2}(\pmb{Z})$, the modular group.
\vspace{.2 cm}

$\bullet$ $\pmb{\mathcal{D}} \simeq \pmb{\Gamma}  \backslash \pmb{\mathcal{H}}$, the standard $SL_{2}(\pmb{Z})$ fundamental domain with cusp at $z = i\infty$
\vspace{.2 cm}

$\bullet$ $\pmb{\Gamma}_{\infty} \subset \pmb{\Gamma}$, the subgroup of upper triangular matrices.
\vspace{.2 cm}

$\bullet$ $\pmb{\zeta}(s)  =   \sum_{n \in \mathbb{N}} n^{-s}$, $\Re(s) > 1$, the Riemann zeta function.
\vspace{.2 cm}

$\bullet$ $\pmb{\zeta}^{*}(s) = \pi^{-s/2}\pmb{\Gamma}(s/2)\pmb{\zeta}(s)$.
\vspace{.2 cm}

$\bullet$ $\pmb{\Theta}_{t}(z) = \sum_{(m,n) \in \mathbb{Z}^{2} \backslash \{0 \}   }e^{-\pi t\frac{|mz + n|^{2}}{y}} $.
\vspace{.2 cm}

$\bullet$   $\pmb{E}_{s}(z) =  \frac{1}{2} \sum_{(c,d) \in \mathbb{Z}^{2},  (c,d) = 1} y^{s}|cz + d|^{-2s}$.
\vspace{.2 cm}

$\bullet$   $\pmb{E}^{*}_{s}(z) =  \pi^{-s}\pmb{\Gamma}(s) \sum_{(m,n) \in \mathbb{Z}^{2} \backslash \{0 \}   } y^{s} |mz + n|^{-2s}$.
\vspace{.2 cm}

$\bullet$     $\pmb{\vartheta}_{t}(y)   = \sum_{ n \in \mathbb{N}_{> 0} }e^{-\pi \frac{t}{y} n^{2}}$.
\vspace{.2 cm}

$\bullet$       $\pmb{\mathcal{M}}_{y}[\varphi](s) =  \int_{0}^{\infty} dy y^{s - 1} \varphi(y)$, the Mellin transform of the function $\varphi$.
\vspace{.2 cm}

$\bullet$   $\mathcal{P}[\varphi](z) = \sum_{\gamma \in \pmb{\Gamma}_{\infty}  \backslash \pmb{\Gamma}} \varphi (\Im(\gamma(z)))  $, the Poincar\'e series of the function $\varphi : \mathbb{R}_{> 0} \rightarrow \mathbb{C}$.
\vspace{.2 cm}

$\bullet$    $\pmb{a}_{0}(y) = \int_{0}^{1} dx f(x,y)$, the constant term  of the modular invariant function $f(x,y) = \sum_{n \in \mathbb{Z}}\pmb{a}_{n}(y)e^{2\pi i n x}$.
\vspace{.2 cm}

$\bullet$    $ \langle f, g \rangle_{ \pmb{\Gamma} \backslash \pmb{\mathcal{H}}} = \int_{\pmb{\mathcal{D}}} dx dy y^{-2} \bar{f}(z) g(z)$, the Petersson inner product of the modular invariant functions $f(z)$, $g(z)$.
\vspace{.2 cm}

$\bullet$      $ \langle \varphi, \xi \rangle_{ \pmb{U} \backslash \pmb{\mathcal{H}}}  =  \int_{0}^{\infty}dy y^{-2} \bar{\varphi}(y) \xi(y)$, the inner product
on the space of functions $\pmb{U} \backslash \pmb{\mathcal{H}} \simeq \mathbb{R}_{> 0}   \rightarrow \mathbb{C}  $.
\vspace{1 cm}


Due to $\pmb{\Gamma}_{\infty} = \pmb{U} \cap SL_{2}(\pmb{Z})$ invariance, a modular invariant function $f = f(x,y)$,
can be decomposed in Fourier series in the $x$ variable, $f(x,y) = \sum_{n \in \mathbb{Z}} \pmb{a}_{n}(y) e^{2\pi i n x}$.
The constant Fourier term $\pmb{a}_{0}(y)$ then gives the $f$ average along the horocycle $\pmb{\mathcal{H}}_y : = (\mathbb{R} + i y)/\pmb{\Gamma}_{\infty}$

\vspace{.27 cm}

\beq
\pmb{a}_{0}(y) = \int_{0}^{1}dx f(x,y) =  \frac{1}{L(\pmb{\mathcal{H}}_y )} \int_{\pmb{\mathcal{H}}_y } d s f,  \label{constantmap}
\eeq

\vspace{.27 cm}

where  $L(\pmb{\mathcal{H}}_y ) = 1/y$ is the horocycle  length, measured by the hyperbolic $\pmb{\mathcal{H}}$ metric, $ds = y^{-1}\sqrt{dx^2 + dy^2 }$.

\vspace{.35 cm}

In this paper  we  focus on two classes of growing conditions for  $SL_{2}(\pmb{Z})$-invariant  functions:

modular functions with polynomial growth at the cusp  $y \rightarrow \infty$
\beq
\pmb{\mathcal{C}}_{TypeII} =\{ f(x,y) \sim \sum_{i = 1}^{l}\frac{c_{i}}{n_{i}!}\, y^{\alpha_{i}} \log^{n_i}y, \, \, y \rightarrow \infty, \,\, \, c_i, \alpha_i \in \mathbb{C}, \Re(\alpha_i )< 1/2,   n_i \in \mathbb{N}_{\ge 0}    \}, \label{poly}
\eeq

and   modular functions with bounded exponential growth at the cusp, whose
  horocycle average $\pmb{a}_{0}(y)$ grows polynomially   at the cusp $y \rightarrow \infty$:

\vspace{.2 cm}

\beq
 \pmb{\mathcal{C}}_{Heterotic} = \{  f(x,y) \sim y^{\alpha}e^{\pi \beta y}e^{2\pi i \kappa x}, y \rightarrow \infty;  \beta < 1, \kappa \in \mathbb{Z} \backslash \{ 0 \}, \Re(\alpha) < 1/2 \}. \label{cstring}
\eeq

\vspace{.2 cm}

The choices of symbols  $\pmb{\mathcal{C}}_{TypeII}$ and  $\pmb{\mathcal{C}}_{Heterotic}$,
reflect   the appearance of modular functions with such growing conditions respectively  in type II
string  and heterotic string  genus one closed string
amplitudes,  (with no tachyons in the spectrum).  Bounds on $\alpha$ and on $\beta$ in (\ref{poly}) and (\ref{cstring})
are universal in string theory, and follow by
consistency requirements,  (unitarity  of the quantum worldsheet conformal field theory [GSW]).

\vspace{.35 cm}

   String theory suggests that  automorphic functions with  growing conditions in
   $\pmb{\mathcal{C}}_{TypeII}$  or in  $\pmb{\mathcal{C}}_{Heterotic}$  do have  convergent  horocycle average
   in the long limit $y \rightarrow 0$, and should exhibit asymptotic behavior  similar to (\ref{theta2})\footnote{
Those hints  follow from the following considerations:
the exponentially growing part for a modular function $f$ in $\pmb{\mathcal{C}}_{Heterotic}$ in string theory language  corresponds to a ''non-physical tachyon '',
a tachyonic state which is not  in the physical spectrum. Indeed,
the exponentially growing part $f(x,y) \sim e^{2\pi i \kappa x} e^{2\pi \beta y}$, $y \rightarrow \infty$,  $\kappa \in \mathbb{Z} \setminus \{ 0 \}$
does not contribute to the $f$   horocycle average, since
\beq
\int_{0}^{1}dx \,  e^{2\pi i \kappa x} e^{2\pi \beta y} = 0, \qquad \kappa \in \mathbb{Z} \setminus \{ 0 \}. \nn
\eeq
Non-physical tachyonic states are expected not to influence the closed string  physical properties. Therefore, one expects
both Type II and Heterotic  strings to have the same  qualitative asymptotic behavior of the spectrum,
i.e. both to  enjoy asymptotic supersymmetry in the absence of \emph{physical} tachyons in their spectra [KS].
This translates back in the expectation for modular functions in   both $\pmb{\mathcal{C}}_{TypeII}$ and $\pmb{\mathcal{C}}_{Heterotic}$
to have the  same asymptotic for their long horocycle average $\pmb{a}_{0}(y)$ in the $y \rightarrow 0$ limit.}.

\vspace{.35 cm}

In this paper  we prove theorems  for  long horocycle average asymptotic of automorphic functions in  $\pmb{\mathcal{C}}_{TypeII}$.
 We also  prove some weaker results for $\pmb{\mathcal{C}}_{Heterotic}$,
and  leave open the complete answer   on long  horocycle averages for automorhic functions in $\pmb{\mathcal{C}}_{Heterotic}$.
 We believe this is an interesting open  question, since  peculiar features of the class of function
 $\pmb{\mathcal{C}}_{Heterotic}$ and the bounds  on $\alpha$ and on  $\beta$ for a sufficient condition for convergence of the
 long horocycle average
 do  not  seem to emerge   from the theory
 of automorphic functions.  A complete answer on the horocycle average asymptotic
 for modular function in $\pmb{\mathcal{C}}_{Heterotic}$ would probe the benefit one may actually gain
 by translating the homogenous dynamics horocycle problem in string theory terms.

\vspace{.35 cm}

In the rest of the  introduction, we summarize our results and  illustrate ideas and methods
employed  to derive them. In order to introduce main concepts on which we focus in this paper, we  start
 in the next section  with a brief illustration  on  how the asymptotic displayed  in (\ref{theta2})
for modular function of rapid decay is derived by the Rankin-Selberg method [Za1],
(more material on that is presented in \S \ref{SPrapiddecay}).

We then switch to modular functions of polynomial growth and discuss why their long horocycle asymptotic behavior
cannot be derived by the \emph{standard} Rankin-Selberg method. In dealing with modular functions of polynomial growth,  Zagier [Za2]
has designed a   Rankin-Selberg
method which is based on an  unfolding method for modular integrals on a truncated version of
the fundamental domain $\pmb{\mathcal{D}}$. We contrast Zagier's method with an alternative unfolding  method we propose here,
which relies on a unfolding trick employing  the theta series  $\pmb{\Theta}_{t}(\tau)$.
 This theta series  $\pmb{\Theta}_{t}(\tau)$ is  related to the spectral Eisenstein series $\pmb{E}_{s}(\tau)$ by a
Mellin transform. One of the  advantages of this Theta  method is to avoid  complications with unfolding  on a truncated version
of the fundamental domain  $\pmb{\mathcal{D}}$.


\vspace{.3 cm}

\subsection{Modular functions of rapid decay and the Rankin-Selberg method}\label{Srapid}
Let us consider the Rankin-Selberg integral

\vspace{.27 cm}

\beq
\langle \pmb{E}_{s}(z), f(z) \rangle_{ \pmb{\Gamma} \backslash \pmb{\mathcal{H}}} = \int_{\pmb{\mathcal{D}}} dx dy y^{-2} \pmb{E}_{s}(z)\, f(x,y) \label{RSint},
\eeq

\vspace{.27 cm}

when    $f =  f(x,y)$ is a modular invariant function  of rapid decay at the cusp $y \rightarrow \infty$.

\vspace{.5 cm}

The spectral Eisenstein series  $\pmb{E}_{s}(z)$ has a Poincar\'e series representation for $\Re(s) > 1$

\vspace{.27 cm}

\beq
\pmb{E}_{s}(z) = \sum_{\gamma \in \pmb{\Gamma}_{\infty} \backslash \pmb{\Gamma}} \Im(\gamma(z))^{s}, \qquad  \Re(s) > 1, \label{partialE}
\eeq

\vspace{.27 cm}

where  $\gamma(z) = \frac{az + b}{cz + d}$, with $\begin{pmatrix} a & b \\ c & d  \end{pmatrix} \in \pmb{\Gamma}$.

The possibility of exchanging the series with the integration on the fundamental domain $\pmb{\mathcal{D}}$

\vspace{.27 cm}

\eqn
 \int_{\pmb{\mathcal{D}}} dx dy y^{-2} \pmb{E}_{s}(z)\, f(x,y) &=& \int_{\pmb{\mathcal{D}}} dx dy y^{-2} \, f(x,y) \, \sum_{\pmb{\gamma} \in \pmb{\Gamma}_{\infty} \backslash \pmb{\Gamma}} \Im(\gamma(z))^{s} \nn \\
&=& \sum_{\gamma \in \pmb{\Gamma}_{\infty} \backslash \pmb{\Gamma}} \int_{\pmb{\mathcal{D}}} dx dy y^{-2} \, f(x,y) \, \Im(\gamma(z))^{s}, \nn \\
\label{unfoldrapid}
\feqn

\vspace{.27 cm}

amounts in being able to perform the unfolding trick. This corresponds to  using modular transformations $\gamma \in  \pmb{\Gamma}_{\infty} \backslash \pmb{\Gamma} $, to unfold the integration domain $\pmb{\mathcal{D}} \simeq \pmb{\Gamma} \backslash \pmb{\mathcal{H}}$ into
the half-infinite strip  $\pmb{\Gamma}_{\infty} \backslash \pmb{\mathcal{H}} \simeq   [-1/2, 1/2 ) \times (0, \infty) \subset \pmb{\mathcal{H}}$.

\vspace{.5 cm}

When $f = f(x,y)$  is of rapid decay at the cusp $y \rightarrow \infty$,  since  $\pmb{E}_{s}(z)$ is of polynomial growth
at the cusp

\vspace{.27 cm}

\beq
\pmb{E}_{s}(z) \sim  y^s + \frac{\pmb{\zeta}^{*}(2s - 1)}{\pmb{\zeta}^{*}(2s)}y^{1 - s} + o(y^{-N}),  \qquad  y \rightarrow \infty, \qquad \forall N > 0, \nn
\eeq

\vspace{.27 cm}

  eq. (\ref{unfoldrapid}) follows   by    Lebesgue dominated convergence theorem on the sequence of products of partial sums of the series in (\ref{partialE})
times the function $f(x,y)$.

This leads to connect the Rankin-Selberg integral to  the Mellin  transform of the function $\pmb{a}_{0}(y)/y$
\beq
\int_{0}^{\infty}dy \, y^{s - 2} \pmb{a}_{0}(y) = \int_{\pmb{\mathcal{D}}} d x dy \, y^{-2}  \pmb{E}_{s}(z)f(x,y). \label{unfl}
\eeq

\vspace{.35 cm}

 A relevant issue at this point  is to determine  analytic properties  of the integral in the  r.h.s. as a function of the complex variable $s$.
 Uniform convergence for $y \rightarrow \infty$ of the  Rankin-Selberg integral
with respect to the complex variable $s$,  assures  that the integral function in the r.h.s. $\pmb{I}(s) : = \langle \pmb{E}_{s}(z), f(z) \rangle_{ \pmb{\Gamma} \backslash \pmb{\mathcal{H}}}$  inherits analytic properties of $\pmb{E}_{s}(z)$.
In the present case  $f$ is of rapid decay, and  uniform convergence of the integral function $\pmb{I}(s)$
holds. Thus the Mellin transform in the l.h.s. of (\ref{unfl})
 inherits as a function of the variable $s \in \mathbb{C}$ the same analytic properties of the Eisenstein series $\pmb{E}_{s}(z)$.

\vspace{.35 cm}

The spectral Eisenstein series  $\pmb{E}_{s}(z)$  has a simple pole in  $s = 1$ with residue $\frac{1}{2\pmb{\zeta}^{*}(2)} = \frac{3}{\pi}$, and poles in $s = \frac{\rho}{2}$,
where $\rho$'s are the non trivial zeros of the Riemann zeta function.

This leads to the following  meromorphic continuation  for the Mellin transform of the function $\pmb{a}_{0}(y)/y$

\beq
\langle y^{s}, \pmb{a}_{0}(y) \rangle_{U \backslash \pmb{\mathcal{H}} } = \int_{0}^{\infty}dy \, y^{s - 2} \pmb{a}_{0}(y) = \frac{C_{0}}{s - 1} + \sum_{\pmb{\zeta}^{*}(\rho) = 0} \frac{C_{\rho}}{s - \rho/2}, \label{RSpoles}
\eeq
where
\beq
C_{0} = Res_{s \rightarrow 1} \int_{\pmb{\mathcal{D}}} dx dy\, y^{-2} \pmb{E}_{s}(z) f(z) = \frac{3}{\pi} \int_{\pmb{\mathcal{D}}} dx dy\, y^{-2} f(z), \nn
\eeq
and
\beq
C_{\rho} = Res_{s \rightarrow \rho/2} \int_{\pmb{\mathcal{D}}} dx dy\, y^{-2} \pmb{E}_{s}(z) f(z), \nn
\eeq

(whenever $\rho$ is a multiple non trivial zero of $\zeta(s)$, one has to raise the denominator in (\ref{RSpoles}) to a power  equal  the order of this zero ).

\vspace{.5 cm}

One finally  obtains the $y \rightarrow 0$ behavior
of $\pmb{a}_{0}(y)$ displayed  in (\ref{theta2})  by using  the  meromorphic continuation  given in  (\ref{RSpoles}),
 whenever the inverse Mellin transform exists, with the help of the following proposition:

\vspace{.27 cm}

\begin{propo}\label{propomellin1}
Let $\varphi = \varphi(y)$  be a function  $\varphi :  (0, \infty) \rightarrow \mathbb{C}$, of rapid decay for $y \rightarrow \infty$,
with Mellin transform $\pmb{\mathcal{M}}[\varphi](s)$.

Suppose, that $\pmb{\mathcal{M}}[\varphi](s)$ can be analytically continued
 to the  meromorphic function
\beq
\pmb{\mathcal{M}}[\varphi](s) = - \sum_{i=1}^{l} \frac{1}{(\alpha_{i} - s)^{n_{i} + 1}}, \qquad \alpha_{i} \in \mathbb{C}, \quad  n_{i} \in \mathbb{N}_{\ge 0},
\eeq
then the following asymptotic holds true
$$ \varphi(y) \sim \sum_{i=1}^{l} \frac{1}{n_{i} !} \, y^{- \alpha_{i} }\log^{n_{i}} y  + o(y^N)   \qquad    y \rightarrow 0, \qquad \forall N > 0. $$
\end{propo}

\vspace{.27 cm}

Therefore, if one supplies extra conditions
 on $f$, which guarantee convergence of the inverse Mellin transform integral,
  (discussion on this matter is postponed to  section \S \ref{SPrapiddecay}),
 then  from eq.  (\ref{RSpoles}) and proposition \ref{propomellin1}, one can prove the  asymptotic  (\ref{theta2}) to hold.

  In section \S \ref{SPrapiddecay}  extra material on the rapid decay case is provided. There,  we also contrast  horocycle average asymptotic
  of $f$ of rapid decay with asymptotic and error estimate  of the rate
  of uniform distribution of the horocycle itself  $\pmb{\Gamma}_{\infty} \backslash  (\mathbb{R} + iy )$    in $\pmb{\mathcal{D}}$
  in the limit $y \rightarrow 0$.

  \vspace{.35 cm}

\subsection{Modular functions of  not-so-mild growing conditions}
  Let us start by discussing what does not go through in the  analysis presented in the previous section when one considers
  modular functions which decay slower at the cusp  then those of rapid decay.

\vspace{.5 cm}

When $f$ is  in $\pmb{\mathcal{C}}_{Type II}$ (\ref{poly}),
the Rankin-Selberg integral in (\ref{RSint}) is convergent for $Max \{\alpha_i \} < \Re(s) < 1 - Max\{\alpha_i \}$,
but it is not uniformly convergent.  When ${\rm{min}}\{ \alpha_{i} \} > 0$ this domain of convergence  is disjointed
from  the strip $\Re(s) > 1$ of convergence of  $\pmb{E}_{s}(z)$  as the  Poincar\'e series (\ref{partialE}).
This implies that one cannot use Lebesgue   dominate convergence theorem
for proving the unfolding trick  (\ref{unfoldrapid}), and thus one cannot reach eq. (\ref{unfl}).

Moreover, for $f \in \pmb{\mathcal{C}}_{Type II}$ , the Rankin-Selberg integral is not uniformly
convergent for $y \rightarrow \infty$ with respect to the complex parameter $s$.
This leads to the expectation that $\pmb{I}(s)$ does not inherits \emph{only} analytic  properties of
$\pmb{E}_{s}(z)$, but that $\pmb{I}(s)$ had singularities also depending on $\alpha_i$, $n_i$.

\vspace{.5 cm}

 Zagier  [Za2] has designed  a Rankin-Selberg method for automorphic functions of polynomial behavior at the cusp
  by devising an unfolding trick for modular integral  restricted to a   truncated version of the fundamental domain
 $\pmb{\mathcal{D}}_T : = \{ x + iy \in \pmb{\mathcal{D}} | y \le T, \,  T > 1 \}$.
 In this way, He  connects   analytic properties of the Rankin-Selberg integral on $\pmb{\mathcal{D}}_T$,
 to various quantities involving the  modular function $f(x,y)$, and its  constant term $\pmb{a}_{0}(y)$.
 Then  by  studying the  $T \rightarrow \infty$ limit, He obtains   analytic properties of
the following  Rankin-Selberg  integral transform

\vspace{.27 cm}

\beq
\pmb{R}^{*}(f,s) := \pmb{\zeta}^{*}(2s) \int_{0}^{\infty}dy y^{s - 2}\left(\pmb{a}_{0}(y) - \varphi(y)\right), \label{RS}
\eeq

\vspace{.27 cm}

where,

\vspace{.27 cm}

$$ \varphi(y) : = \sum_{i = 1}^{l}\frac{c_{i}}{n_{i}!}\, y^{\alpha_{i}} \log^{n_i}y$$ is the leading  polynomial growing  part of $f(x,y)$
  in the $y \rightarrow \infty$ limit.

\vspace{.27 cm}

 $\pmb{R}^{*}(f,s)$ is the relevant integral transform for the polynomial growth case, which parallels
 the Mellin transform (\ref{unfl}) of the rapid decay case.
Analytic continuation of $\pmb{R}^{*}(f,s)$ is given by the following theorem:

\vspace{.3 cm}

\begin{theorem}($\rm \textbf{Zagier}$, \cite{Za2})\label{Zpol}
Let $f$ be  a modular invariant function of  polynomial growth at the cusp
$$f(x,y) \sim \sum_{i = 1}^{l}\frac{c_{i}}{n_{i}!}\, y^{\alpha_{i}} \log^{n_i}y + o(y^{- N}), \qquad y \rightarrow \infty,  \qquad \forall N > 0,$$
then  the Rankin-Selberg transform (\ref{RS}) can be analytically continued to the  meromorphic function

\beq
\pmb{R}^{*}(f,s) = \sum_{i=1}^l c_i \left( \frac{\pmb{\zeta}^{*}(2s)}{(1 - s - \alpha_i)^{n_i + 1}}
 + \frac{\pmb{\zeta}^{*}(2s - 1)}{(s - \alpha_i )^{n_i + 1}}  + \frac{{\rm{ entire \,\, function \,\, of} \,\, s}}{s(s - 1)}    \right). \label{eq-RSZpolar}
\eeq
\end{theorem}

\vspace{.3 cm}

Eq. (\ref{eq-RSZpolar})  parallels eq. (\ref{RSpoles})
of the rapid decay case.

\vspace{.5 cm}

 We shall now  present our methods, which allow also to prove theorem \ref{Zpol} by a distinct route. This route  avoids to
use unfolding tricks on truncated versions of $\pmb{\mathcal{D}}$ as in [Za2].
With this method  we will also prove various  results of this paper.
In order to illustrate our methods, and in the polynomial growing case, to contrast it with those in
[Za2], we start by introducing   the following


{\centering
\subsection*{\emph{ Lattices series magic square}:}

\vspace{.5 cm}
\beq
\begin{array}[c]{ccc}
\pmb{\Theta}_{t}(z) &      \overset{\pmb{\mathcal{M}}_{t} }{\longrightarrow }&  \pmb{E}^{*}_{s}(z)    \\
  &  &  \\
\uparrow \mathcal{P}_y &                  &   \uparrow  \mathcal{P}_{y}  \\
& & \\
 \pmb{\vartheta}_{t}(\Im(z)) &   \overset{\pmb{\mathcal{M}}_{t} }{\longrightarrow } &  \pmb{\mathcal{E}}^{*}_{s}(z)      
\end{array} \label{magico}
\eeq

\vspace{.5 cm}}

relating    four  functions of great relevance  in analytic number theory.
In the upper vertexes of the  square sit two $2$-dimensional lattices series,
the dressed spectral Eisenstein series  $\pmb{E}^{*}_{s}(z)$,
and the $2$-lattice theta series $\pmb{\Theta}_{t}(z)$,

\beq
\pmb{E}^{*}_{s}(z) : = \pi^{-s} \pmb{\Gamma}(s) \sum_{\omega \in \Lambda_z} \left( \frac{ |\omega|^{2}}{\Im(z)} \right)^{-s}, \nn
\eeq

\beq
\pmb{\Theta}_{t}(z) : = \sum_{\omega \in \Lambda_z} e^{ - \pi t \left(\frac{ |\omega|^{2}}{\Im(z)}\right)}, \nn
\eeq

with $\Lambda_z : = \{ m z + n \in \mathbb{C}, \, (m,n) \in \mathbb{Z}^2  \setminus \{ 0 \}, z \in \pmb{\mathcal{H}}  \}$ a two
dimensional lattice, with modular parameter $z$.
These two $2$-lattice series are related  by Mellin integral transform $\pmb{\mathcal{M}}$

\beq
\pmb{E}^{*}_{s}(z) : = \int_{0}^{\infty} dt \,  t^{s-1} \pmb{\Theta}_{t}(z)   =  \int_{0}^{\infty} dt \,  t^{s-1}\sum_{\omega \in \Lambda_z}e^{ - \pi t \left(\frac{ |\omega|^{2}}{\Im(z)}\right)}. \nn
\eeq

In the lower vertices of the magic square sit two  $1$-dimensional lattice  $\mathbb{N}_{> 0}$ series,
that are the homologous of the two dimensional ones

\beq
\pmb{\mathcal{E}}^{*}_{s}( \Im(z) ) : =  \pi^{-s} \pmb{\Gamma}(s) \sum_{n \in  \mathbb{N}_{> 0}} \left( \frac{ n^{2}}{\Im(z)} \right)^{- s} = \Im(z)^s \pmb{\zeta}^{*}(2s),  \nn
\eeq

\beq
\pmb{\vartheta}_{t}(\Im(z)) : =  \sum_{n \in  \mathbb{N}_{> 0}}e^{ - \pi t \left(\frac{ n^{2}}{\Im(z)}\right)}. \nn
\eeq

The above  two $1$-dimensional lattice series are also related by a Mellin integral transform

\beq
\pmb{\mathcal{E}}^{*}_{s}(\Im(z)) : = \int_{0}^{\infty} dt \,  t^{s-1} \pmb{\vartheta}_{t}(\Im(z)).
\eeq

The vertical arrows in the magic square uplift one dimensional lattice series to two dimensional lattice series.
This works through the relation $\Lambda_z = \mathbb{N}_{> 0} \otimes  \tilde{\Lambda}_{z}$, where
$\tilde{\Lambda}_{z} : = \{ cz + d | (c,d) \in  \mathbb{Z}^2,  (c,d) = 1 \}$ is the co-primed  $2$-lattice.
  The $\tilde{\Lambda}_z$ modular group is   $\pmb{\Gamma} \sim SL_{2}(\pmb{Z})$ identified by the $\mathbb{N}_{> 0}$ left action, i.e. $\pmb{\Gamma}_{\infty} \backslash \pmb{\Gamma}$.
Therefore
\eqn
\pmb{E}^{*}_{s}(z)  &=& \pi^{-s} \pmb{\Gamma}(s) \sum_{\omega \in \Lambda_z} \left( \frac{ |\omega|^{2}}{\Im(z)} \right)^{-s}, \nn \\
                     &=& \pi^{-s} \pmb{\Gamma}(s)\sum_{n \in \mathbb{N}_{> 0}} \sum_{\tilde{\omega} \in \tilde{\Lambda}_z} \left( \frac{ |\tilde{\omega}|^{2}}{\Im(z)} \right)^{-s} \nn \\
                     &=&   \sum_{\gamma \in \pmb{\Gamma}_{\infty} \backslash \pmb{\Gamma}} \mathcal{E}^{*}_{s}(\Im(\gamma(z))), \nn
\feqn

and by applying a reasoning as above

\beq
\pmb{\Theta}^{*}_{t}(z) =  \sum_{\gamma \in \pmb{\Gamma}_{\infty} \backslash \pmb{\Gamma}} \pmb{\vartheta}_{t}(\Im(\gamma(z))).
\eeq

\vspace{.3 cm}

Given a modular invariant function $f = f(x,y)$,  by taking inner products both in  $ \pmb{\Gamma} \backslash \pmb{\mathcal{H}}$
aand in  $ \pmb{U} \backslash \pmb{\mathcal{H}} \simeq \mathbb{R}_{> 0}$ with  functions appearing  in  diagram (\ref{magico}), one finds a set of relations displayed  by  the following

{\centering
 \subsection*{\emph{Inner products magic square}:}

\vspace{.5 cm}

\beq
\begin{array}[c]{ccc}
\langle \pmb{\Theta}_{t}(z), f(z) \rangle_{ \pmb{\Gamma} \backslash \pmb{\mathcal{H}}}    &      \overset{\pmb{\mathcal{M}}_{t} }{\longrightarrow }& \langle \pmb{E}^{*}_{s}(z), f(z) \rangle_{ \pmb{\Gamma} \backslash \pmb{\mathcal{H}}}   \\
& & \\
\downarrow Unfolding  &                  &   \downarrow   Unfolding  \\
&  & \\
\langle \pmb{\vartheta}_{t}(y), \pmb{a}_{0}(y) \rangle_{U \backslash \pmb{\mathcal{H}} } &   \overset{\pmb{\mathcal{M}}_{t} }{\longrightarrow } &  \pmb{\zeta}^{*}(2s)\langle y^{s}, \pmb{a}_{0}(y) \rangle_{U \backslash \pmb{\mathcal{H}} }.
\end{array} \label{magicp}
\eeq

\vspace{.5 cm}

}

The inner product on $ \pmb{\Gamma} \backslash \pmb{\mathcal{H}} $ corresponds to the Petersson inner product between two modular invariant functions.
 Given two modular functions  $f$ and $g$, it  is defined as follows
\beq
\langle f, g \rangle_{ \pmb{\Gamma} \backslash \pmb{\mathcal{H}}} : = \int_{ \pmb{\Gamma} \backslash \pmb{\mathcal{H}}} dx dy \,y^{-2} \bar{f}(z) g(z), \label{Gprod}
\eeq
where $\bar{f}$ is the complex conjugate of $f$.
The inner product on  $\pmb{U} \backslash \pmb{\mathcal{H}}$ for
 a pair of  functions $\varphi$ and $\xi$ on $\mathbb{R}_{>0}$ with values in $\mathbb{C}$ is defined as

\beq
\langle \varphi, \xi \rangle_{\pmb{U} \backslash \pmb{\mathcal{H}}} : = \int_{0}^{\infty} dy \, y^{-2} \bar{\varphi}(y) \xi(y). \label{Uprod}
\eeq

\vspace{.3 cm}

Vertical arrows in the diagram (\ref{magicp}) correspond to the following  unfolding trick, which
 allows  to identify  the constant map $\pmb{a}_{0}$  as the   adjoint map of the  Poincar\'e map  $\mathcal{P}$  with respect to the
  inner products  (\ref{Gprod}) and   (\ref{Uprod})

\eqn
 \langle \mathcal{P}[\varphi], f \rangle_{ \pmb{\Gamma} \backslash \pmb{\mathcal{H}}}
&=& \int_{ \pmb{\Gamma} \backslash \pmb{\mathcal{H}}} dx dy \,y^{-2} \mathcal{P}[\varphi](z) f(z), \nn \\
&=& \int_{ \pmb{\Gamma} \backslash \pmb{\mathcal{H}}} dx dy \,y^{-2} f(z) \, \sum_{\gamma \in  \pmb{\Gamma}_{\infty} \backslash \pmb{\Gamma}} \varphi(\Im(\gamma(z))) \nn \\
&=& \sum_{\gamma \in  \pmb{\Gamma}_{\infty} \backslash \pmb{\Gamma}} \int_{ \Gamma \backslash \pmb{\mathcal{H}}} dx dy \,y^{-2} f(z) \, \varphi(\Im(\gamma(z))) \nn \\
&=& \int_{0}^{\infty} dy \, y^{-2} \bar{\varphi}(y) \int_{0}^{1} dx f(z)  \nn \\
&=& \langle  \varphi ,  \pmb{a}_{0}[f]  \rangle_{U \backslash \pmb{\mathcal{H}}}, \nn \\ \label{utrick}
\feqn

where $\pmb{a}_{0}[f]$ is the constant map,
\beq
\pmb{a}_{0}[f](y) : = \int_{0}^{1} dx f(x,y). \nn
\eeq
 As already remarked in (\ref{constantmap}),  the constant map $\pmb{a}_{0}$ in geometrical terms  gives   the horocycle average of the modular invariant function $f$.

 The above unfolding  trick is equivalent of being able to
exchange in the inner product $\langle \mathcal{P}[\varphi], f \rangle_{ \pmb{\Gamma} \backslash \pmb{\mathcal{H}}}$  the series over modular transformations in $\pmb{\Gamma}_{\infty} \backslash \pmb{\Gamma}$, with   integration
on the fundamental domain $\pmb{\mathcal{D}} \simeq \pmb{\Gamma}  \backslash \pmb{\mathcal{H}}$.
This possibility  depends   on the behavior at the cusp  of the product of  the modular function $f(z)$  with
the Poincar\'e series  $\sum_{\gamma \in  \pmb{\Gamma}_{\infty} \backslash \pmb{\Gamma}} \varphi(\Im(\gamma(z)))$.

\vspace{.5 cm}

In the rest of this introduction, we discuss and contrast the classical
Rankin-Selberg method, which we introduced in  \S \ref{Srapid}, and it  corresponds to moving along the \emph{right} column of diagram
\ref{magicp} in the direction of the arrow, to a Theta unfolding method. This alternative method   corresponds to moving along the \emph{left}
column of diagram \ref{magicp} in the direction of the vertical arrow, and then by using the horizontal lower arrow.
For various  classes of growing conditions at the cusp,
we shall contrast unfolding of a modular integral of the product  of a  function $f$  with the spectral  Eisenstein series $\textbf{E}^{*}_{s}(z)$,
with unfolding by using the  double theta series $\pmb{\Theta}_{t}(z)$.
 Discussions and results of this paper should illustrate  advantages of using  the double theta series  $\pmb{\Theta}_{t}(z)$
 unfolding trick, when one considers modular invariant functions which have not-so-mild growing conditions at the cusp.
The general idea is that whether $\pmb{E}^{*}_{s}(z)$ grows polynomially at the cusp, $\pmb{\Theta}_{t}(z)$ provides a better
  convergence for the modular integral, since the subseries of terms of $\pmb{\Theta}_{t}(z)$
which decay exponentially at the cusp, are precisely those   which allow to perform the unfolding trick.
This  unfolding trick allows a better control for  modular functions with not-so-mild growing condition at the cusp.


To summarize,
our Theta method corresponds to  the following route in  the \ref{magicp} diagram
\vspace{.5 cm}

\beq
\begin{array}[c]{ccc}
\langle \pmb{\Theta}_{t}(z), f(z) \rangle_{ \pmb{\Gamma} \backslash \pmb{\mathcal{H}}}    &    &    \\
& & \\
\downarrow Unfolding  &                     &  \\
&  & \\
\langle \pmb{\vartheta}_{t}(y), \pmb{a}_{0}(y) \rangle_{U \backslash \pmb{\mathcal{H}} } &   \overset{\pmb{\mathcal{M}}_{t} }{\longrightarrow } &  \pmb{\zeta}^{*}(2s)\langle y^{s}, \pmb{a}_{0}(y) \rangle_{U \backslash \pmb{\mathcal{H}} }  = \pmb{R}^{*}(f,s).
\end{array} \label{magic2}
\eeq

\vspace{.5 cm}
For $f$ of polynomial growth, the advantage of this route is that no  truncations of the domain of integration $\pmb{\mathcal{D}}$
are required.
Unfolding of the integration domain in the modular integral

\vspace{.27 cm}

\beq
\langle \pmb{\Theta}_{t}(z), f(z) \rangle_{ \pmb{\Gamma} \backslash \pmb{\mathcal{H}}} = \int_{\pmb{\mathcal{D}}} dx dy \, y^{-2} f(x,y) \, \sum_{(m,n) \in \mathbb{Z}^2 \setminus \{ 0 \}} e^{- \frac{\pi}{y}|m z +  n|^{2}}, \label{TheraIntP}
\eeq

\vspace{.27 cm}

 follows from the following  decomposition for the  theta series $\pmb{\Theta}_{t}(z)$

 \beq
 \pmb{\Theta}_{t}(z) = \sum_{(m,n) \in \mathbb{Z}^2 \setminus \{ 0 \}}e^{-\pi t \frac{(mx + n)^{2} + m^2 y^2 }{y}} =
  1 + \pmb{\vartheta}_{t}(\Im(z)) + \sum_{\gamma \in \pmb{\Gamma}_{\infty} \backslash  \pmb{\Gamma}^{'}} \pmb{\vartheta}_{t}(\Im(\gamma(z))),
  \qquad \pmb{\Gamma}^{'} : =  \pmb{\Gamma} \, \, \setminus \, \, \{ \mathbb{I} \},  \label{Theta}
 \eeq
with $\gamma(z) : = \frac{az + b}{cz + d}$.

\vspace{.35 cm}

 One uses modular transformations appearing in  the third term on the r.h.s.,
 of   the form

\beq
 \pmb{\vartheta}_{t}(\Im(\gamma(z))) = \sum_{r \ne 0}e^{-\pi \frac{r^2 t}{y} ((c x + d)^{2} + c^2 y^2)} \qquad c, d \in \mathbb{Z}, c \ne 0, (c,d) = 1, \nn
 \eeq

that  correspond  to the $m \ne 0$ subseries in (\ref{Theta}), and  decay exponentially for $y \rightarrow \infty$.
 Thus, for $f$ in $\pmb{\mathcal{C}}_{TypeII}$, by dominate convergence theorem one can
unfold the modular integral $\langle \pmb{\Theta}_{t}(z), f(z) \rangle_{ \pmb{\Gamma} \backslash \pmb{\mathcal{H}}} $
in the upper vertex of the triangular diagram  \ref{magic2},
and obtain the quantity in the left lower vertex  $\langle \pmb{\vartheta}_{t}(y), \pmb{a}_{0}(y) \rangle_{U \backslash \pmb{\mathcal{H}} }$.
This corresponds to prove  the  vertical arrow  of the triangular  diagram \ref{magic2}  to hold for functions in
 $\pmb{\mathcal{C}}_{TypeII}$.

\vspace{.35 cm}

As a next step, in section \ref{SPpolynomial}, we estimate both
the $t \rightarrow 0$ and the  $t \rightarrow \infty$ asymptotics of the function

\vspace{.27 cm}

\beq
\pmb{i}(t) : = \langle \pmb{\vartheta}_{t}(y), \pmb{a}_{0}(y) \rangle_{U \backslash \pmb{\mathcal{H}} }, \label{iFunction}
\eeq

\vspace{.27 cm}

which appears in the left lower vertex of \ref{magic2}. Due to the  arrow in the lower side of the triangular diagram \ref{magic2},
knowledge  of  $t \rightarrow 0$ and $t \rightarrow \infty$ asymptotics  of the function $\pmb{i}(t)$ allows to reconstruct
meromorphic  expansion of its Mellin transform in  the right lower vertex of \ref{magic2}.
Since the function in the right lower corner  coincides  with the Rankin-Selberg transform of
the constant term $\pmb{a}_{0}(y)$, this allows  to prove theorem \ref{Zpol}.
 Moreover,  the lower row of diagram  (\ref{magic2})
shows  a simple connection between the two functions $\pmb{i}(t)$ and $\pmb{a}_{0}(y)$. This allows to obtain  the $y \rightarrow 0$
asymptotic of $\pmb{a}_{0}(y)$ by means of  the $\pmb{i}(t)$ asymptotic. By this route, in \S \ref{SPpolynomial} we shall  prove the following

\vspace{.5 cm}

\begin{theorem} \label{thpoly2}
For a  given  $f = f(x,y)$  modular invariant function  with polynomial behavior at the cusp
$$
f(x,y) \sim  \sum_{i = 1}^{l}\frac{c_{i}}{n_{i}!}\, y^{\alpha_{i}} \log^{n_i}y  +  o(y^{-N}), \qquad y \rightarrow \infty  \qquad   \forall N > 0,
$$

for $c_i, \alpha_i \in \mathbb{C}$, $\Re(\alpha_i ) < 1/2$, $n_i \in \mathbb{N}_{\ge 0} $, the following asymptotic holds true

$$ \pmb{a}_{0}(y) \sim C_{0} + \sum_{\pmb{\zeta}^{*}(\rho) = 0} C_{\rho} y^{1 - \frac{\rho}{2}} +  \sum_{i = 1}^{l}\frac{c_{i}}{n_{i}!}\,  \frac{\pmb{\zeta}^{*}(2\alpha_{i} - 1)}{\pmb{\zeta}^{*}(2\alpha_{i} )}y^{1 - \alpha_i }\log^{n_i} y + o(y^{N}), \quad y \rightarrow 0, \quad \forall N > 0, $$

where
$$C_{0} = \frac{3}{\pi}\int_{\pmb{\mathcal{D}}}dx dy y^{-2} f(z).$$
\end{theorem}

\vspace{.5 cm}

We now sketch how in \S \ref{SPpolynomial}  we do prove  asymptotics for the function $\pmb{i}(t)$.
This is done in two steps,
first    we need  the following lemma

\vspace{.5 cm}

\begin{lemma}\label{lemma1} Given a  modular invariant function $f = f(x,y)$ with finite integral on $\pmb{\mathcal{D}}$,  $C_0 : = < 1, f >_{ \pmb{\Gamma} \backslash \pmb{\mathcal{H}}}$. Let   $\pmb{a}_{0}(y)$ the $f$ constant Fourier term, then the  following relation holds true

$$\langle \pmb{\vartheta}_{t}(y), \pmb{a}_{0}(y) \rangle_{U \backslash \pmb{\mathcal{H}}} = \frac{1}{t}\langle \pmb{\vartheta}_{1/t}(y), \pmb{a}_{0}(y) \rangle_{U \backslash \pmb{\mathcal{H}}} + \frac{C_{0}}{t} - C_{0}.$$

\end{lemma}

\vspace{.5 cm}

Lemma \ref{lemma1} then allows to  prove the following lemma on  asymptotics of the function  $\pmb{i}(t)$

\vspace{.5 cm}

\begin{lemma}\label{propolemma2}
Let $f =  f(x,y)$ a modular invariant function  with polynomial behavior at the cusp

 $$f(x,y) \sim \sum_{i=1}^{l}\frac{c_i}{n_{i}!} \, y^{\alpha_i}\, \log^{n_i}y  + o(y^{-N}), \qquad y \rightarrow \infty  \qquad \forall N > 0 $$

where $\alpha_i, c_i \in \mathbb{C}$,  $\Re(\alpha_{i}) < 1/2$, $n_i \in \mathbb{N}_{\ge 0}$.

\vspace{.2 cm}

Then, for the function $\pmb{i}(t) : = \langle \pmb{\vartheta}_{t}(y), \pmb{a}_{0}(y) \rangle_{U \backslash \pmb{\mathcal{H}}}$ the following
asymptotics hold true
\eqn
&i)& \qquad \pmb{i}(t) \sim \sum_{i=1}^{l}\frac{c_i}{n_{i}!}  \pmb{\zeta}^{*}(2\alpha_{i} - 1 ) t^{\alpha_{i} - 1}\, \log^{n_i} t + O\left(t^{A  - 1}\log^{N - 1} t   \right),   \qquad  t \rightarrow \infty   \nn \\
&ii)& \qquad  \pmb{i}(t) \sim  - C_{0} +  \frac{C_{0}}{t} - \sum_{i=1}^{l}\frac{c_i}{n_{i}!}    \pmb{\zeta}^{*}(2\alpha_{i} - 1) t^{- \alpha_{i}} \log^{n_i} t
 + O\left(t^{-a}\log^{n - 1} t   \right)  ,  \qquad   t \rightarrow 0 \nn
\feqn
where   $A : = max\{ \Re(\alpha_i) \}$, $a : = min\{ \Re(\alpha_i) \}$, $N : =   max\{ n_i \}$,  $n : =  min\{ n_i \}$, and
$$C_{0} = \int_{\pmb{\mathcal{D}}}dx dy y^{-2} f(z).   $$

\end{lemma}

\vspace{.5 cm}

Then we also recover  Zagier  result on the analytic continuation   (\ref{eq-RSZpolar}) of the Rankin-Selberg transform $\pmb{R}^{*}(f,s)$ in theorem \ref{Zpol}.
Our proof for (\ref{eq-RSZpolar}) uses  lemma \ref{propolemma2},  the horizontal arrow   in diagram \ref{magic2}  and proposition \ref{propomellin1}.
Thereafter, with all the collected results we prove  theorem \ref{thpoly2},
 on  the long horocycle average asymptotic of functions in $\pmb{\mathcal{C}}_{TypeII}$.

\vspace{.3 cm}

\subsection{String inspired class of modular functions of exponential growth  at the cusp}\label{Sexp}

 Section  \S \ref{SPexponential} deals with   the  class of modular function with (bounded) exponential  growing conditions in  $\pmb{\mathcal{C}}_{Heterotic}$ (\ref{cstring}). Examples of functions with such  exponentially growing conditions do
 appear in one-loop amplitudes in heterotic string theory.   We are able to prove much weaker results on the $y \rightarrow 0$ behavior of their horocycle average. However, string theory suggests  better converging  behavior then what we managed to prove in this paper.
 We leave string theory suggestions as  open question at the end of section \S \ref{SPexponential}.
By following the  route  given by  the arrows in  diagram (\ref{magic2}),
we are able to prove the following bound on the growing of the long horocycle average
for modular functions in $\pmb{\mathcal{C}}_{Heterotic}$:

\vspace{.3 cm}

\begin{theorem} \label{thbound}
Let $f = f(x,y)$ a modular invariant function with growing conditions in the class $\pmb{\mathcal{C}}_{Heterotic}$ defined by eq.  (\ref{cstring}), then

\vspace{.27 cm}

\beq
\pmb{a}_{0}(y) \sim o(e^{C/y}) \qquad y \rightarrow 0, \qquad \forall c \in \mathbb{C}, \, \Re(C) > 0.
\eeq
\end{theorem}

\vspace{.27 cm}

  As discussed at the beginning of the introduction, string theory suggests a much stronger result
  on the $\pmb{a}_{0}(y)$ asymptotic, namely  that in the $y \rightarrow 0$ limit  $\pmb{a}_{0}(y)$ is convergent  have  asymptotic as in theorem \ref{thpoly2}.
This leads to the following open question:

\vspace{.5 cm}

\begin{open} \label{Qexp}
(Prove or disprove the following statement):  Given $f = f(x,y)$ modular invariant function in the class $\pmb{\mathcal{C}}_{Heterotic}$ (\ref{cstring}),
 the following asymptotic holds true

\vspace{.27 cm}

 \beq
\pmb{a}_{0}(y) \sim C_{0} + \sum_{\pmb{\zeta}^{*}(\rho) = 0} C_{\rho} y^{1 - \frac{\rho}{2}} + \frac{\pmb{\zeta}^{*}(2\alpha - 1)}{\pmb{\zeta}^{*}(2\alpha)}y^{1 - \alpha} + o(y^{N}), \qquad  \forall N > 0, \qquad y \rightarrow 0, \nn
\eeq

\vspace{.27 cm}

and,

\vspace{.27 cm}

\beq
C_{0} = \frac{3}{\pi}\int_{\pmb{\mathcal{D}}} dy \, y^{-2} \int dx  f(x.y).   \nn
\eeq

\vspace{.27 cm}

where this integral is meant in the conditional sense, with integration along the real axis  performed first.
\end{open}

\vspace{.27 cm}

 Besides discussing string theory  hints for the above open  question,
 at the end of section \ref{SPexponential} we also remark   the possibility of
 having a sort of rigidity in the way the constant term $\pmb{a}_{0}(y)$ may grow in the $y \rightarrow 0$ limit.
 The following  result related to this issue  is given at the end of section \ref{SPexponential}:

\vspace{.3 cm}

\begin{propo} \label{propmod}
Given a $SL_{2}(\pmb{Z})$ invariant function $f$ which grows as $f(x,y) \sim e^{2\pi \beta y}e^{2\pi i \kappa x}$ for $y \rightarrow \infty$
for a certain non-zero  integer  $\kappa \in \mathbb{Z} \,  \backslash \,  \{0 \}$. Then

\vspace{.27 cm}

\beq
\pmb{a}_{0}(y) +  \sum_{r \in \mathbb{Z} \, \backslash \, \{ 0\}}\pmb{a}_{r}(y)e^{2\pi i r \frac{a}{c}} \sim e^{- 2 \pi i \kappa \frac{d}{c}} e^{2\pi \beta \frac{c^{2}}{y}},  \qquad y \rightarrow 0, \nn
\eeq

\vspace{.27 cm}

for every pairs of Farey fractions $\frac{a}{c}$, $\frac{d}{c}$,  $a,c,d \in \mathbb{Z}$, $(a,c) = 1$, $|a| < c$, $(d,c) = 1$
$|d| < c$, $c > 0$.
 $\pmb{a}_{r}(y)$ are the Fourier modes in the expansion $f(x,y) = \sum_{r \in \mathbb{Z}} \pmb{a}_{r}(y)e^{2\pi i r x}$.

\end{propo}

\vspace{.3 cm}

We  end up section \ref{SPexponential} by discussing  the possibility that   proposition \ref{propmod}
together with the bound given by theorem \ref{thbound}
may be of help  in addressing  the open  question raised
in  the open problem \ref{Qexp}.

\vspace{.5 cm}




\section{\Large{\textbf{Rapid decay case: the Rankin-Selberg method and Zagier connection to RH}}}\label{SPrapiddecay}

\vspace{.3 cm}

 This section  contains  a review in some details of the Rankin-Selberg method \newline [R-S] for automorphic functions of rapid decay,
 (some of the material contained in this section overlaps with \S \ref{Srapid}).
 We  review in details,  Zagier proof [Za1] of the
dependence of the error estimate in the horocycle average   asymptotic of modular function of rapid decay on the Riemann hypothesis,
 (eq. (\ref{theta2}) in the introduction).
Most of the material is contained in [Za1], although we have expanded some of the  discussions in [Za1].

\vspace{.3 cm}
Given $f = f(x,y)$ a modular  invariant function of rapid decay at the cusp $y \rightarrow \infty$,
the Rankin-Selberg integral is the following modular integral

\vspace{.27 cm}

\beq
\pmb{I}(s) = \int_{\pmb{\mathcal{D}}}dx dy y^{-2} f(z)  \pmb{E}_{s}(z), \label{modInt}
\eeq

\vspace{.27 cm}

on the $SL_{2}(\pmb{Z})$ fundamental domain $\pmb{\mathcal{D}}$, where

\vspace{.27 cm}

\beq
\pmb{E}_{s}(z)  = \sum_{\gamma \in \pmb{\Gamma}_{\infty} \backslash \pmb{\Gamma}} \Im(\gamma(z))^s =
\frac{1}{2}\sum_{c,d \in  \mathbb{Z}, (c,d)=1}\frac{y^{s}}{|cz + d|^{2s}}, \qquad \Re(s) > 1, \nn
 \eeq

 \vspace{.27 cm}

 is the spectral Eisenstein series. $\pmb{E}_{s}(z)$ can be analytically continued
 to the full plane $s$, except for a simple pole in  $s = 1$ with residue $\frac{3}{\pi}$, and poles
 in $s = \rho/2$, where $\rho$'s are the non trivial zeros of the Riemann zeta function, $\pmb{\zeta}^{*}(\rho) = 0$.

\vspace{.5 cm}

\subsection{Unfolding and analytic heritage}

 The sequence of   partial sums of $\pmb{E}_{s}(z)$ times the function $f(z)$
 is dominated by $\pmb{E}_{s}(z)|f(z)|$, a integrable  function  on  $\pmb{\mathcal{D}}$, for $\Re(s) > 1$.
Thus, by dominated convergence Lebesgue theorem, one can exchange the series with  the integral,
which amounts  to  use the  unfolding trick  for enlarging  the integration domain to
half-infinite  strip  $ [-1/2, 1/2) \times (0, \infty) \subset \pmb{\mathcal{H}}$
\eqn
\pmb{I}(s) &=& \int_{\pmb{\mathcal{D}}}dx dy y^{-2} f(z) \sum_{\gamma \in \pmb{\Gamma}_{\infty} \backslash \pmb{\Gamma}} \Im(\gamma(z))^s \nn \\
&=& \sum_{\gamma \in \pmb{\Gamma}_{\infty} \backslash \pmb{\Gamma}}\int_{\pmb{\Gamma}  \backslash\pmb{\mathcal{H}}}dx dy y^{s - 2} f(z) \nn \\
&=& \int_{0}^{\infty}dy y^{s - 2}\pmb{a}_{0}(y). \nn
\feqn

\vspace{.35 cm}

The integral function $\pmb{I}(s)$ inherits analytic properties of $\pmb{E}_{s}(z)$,
since the modular integral (\ref{modInt}) is uniformly convergent for $y \rightarrow \infty$ in the complex parameter $s$.
 In fact $\pmb{E}_{s}(z)$ grows polynomially for $y \rightarrow \infty$
 \beq
 \pmb{E}_{s}(z) \sim y^s + \frac{\pmb{\zeta}^{*}(2s - 1 )}{\pmb{\zeta}^{*}(2s)} y^{1 - s} \qquad y \rightarrow \infty, \nn
 \eeq
while $f(x,y)$ is of rapid decay for $y \rightarrow \infty$.

\vspace{.35 cm}

Uniform convergence of the modular integral for the complex parameter $s$ on a set $\pmb{A}$ for $z \rightarrow i\infty$  means that given $\epsilon > 0$ there exists a corresponding  neighborhood $\mathcal{U}_{\epsilon}$ of the cusp $z = i\infty$
such that
\beq
\left|\int_{\mathcal{U}_{\epsilon}}dx dy y^{-2} f(z) \partial^{n}_{s}\pmb{E}_{s}(z) \right| < \epsilon \qquad \forall s \in \pmb{A}, \qquad \forall n. \nn
\eeq
 In this case $\mathcal{U}_{\epsilon} = \{z \in \pmb{\mathcal{D}}| \Im(z) >  M_{\epsilon} \}$.

\vspace{.5 cm}

\subsection{\textbf{Poles and Residues of $\pmb{E}_{s}(z)$}}

 $\pmb{E}_{s}(z)$ has a simple  pole in $s = 1$ with residue $3/\pi$.
In fact
\beq
 \pmb{E}^{*}_{s}(z)  = \pmb{\zeta}^{*}(2s)\pmb{E}_{s}(z) =  \frac{1}{2}\pi^{-s}\pmb{\Gamma}(s)\sum_{(m,n) \in  \mathbb{Z}^{2} \backslash \{ 0 \}}\frac{y^{s}}{|m z + n|^{2s}}, \nn
\eeq
is the Mellin transform with respect to the variable $t$ of the function $\pmb{\Theta}_{t}(z)$
\beq
 \pmb{E}^{*}_{s}(z) = \frac{1}{2} \int_{0}^{\infty}dt \, t^{s - 1} \pmb{\Theta}_{t}(z). \label{MelTheta}
\eeq

Double Poisson summation gives

\vspace{.27 cm}

\beq
\pmb{\Theta}_{t}(z) = - 1 + \frac{1}{t}   + \frac{1}{t}\pmb{\Theta}_{1/t}(z), \nn
\eeq
\vspace{.27 cm}

 thus
\vspace{.27 cm}

\beq
\pmb{\Theta}_{t}(z) \sim  - 1  + \frac{1}{t}    \qquad t \rightarrow 0, \nn
\eeq

\vspace{.27 cm}

while $\pmb{\Theta}_{t}(z)$ is of rapid decay for $t \rightarrow \infty$.
Therefore by proposition \ref{propomellin1},  $\pmb{E}^{*}_{s}(z)$ has a pole in $s= 0$ with residue $-1/2$ and
pole in $s = 1$  with residue $1/2$. Therefore,

\vspace{.27 cm}

\beq
\pmb{E}_{s}(z) = \frac{\pmb{E}^{*}_{s}(z)}{\pmb{\zeta}^{*}(2s)} \nn
\eeq

\vspace{.27 cm}

has a pole in $s = 1$ with residue $\frac{1}{2\pmb{\zeta}^{*}(s)} = 3/\pi$ and poles in $\rho/2$, where $\rho$'s
 are the  zeros of the Riemann zeta function $\pmb{\zeta}^{*}(\rho) = 0$.

\vspace{.5 cm}

\subsection{Zagier's result on $\pmb{a}_{0}(y)$,  $y \rightarrow 0$ asymptotic}

 A sufficient condition for the following  $y \rightarrow 0$ asymptotic to hold, (displayed in (\ref{theta2}))

\vspace{.27 cm}

\beq
\pmb{a}_{0}(y) \sim C + \sum_{\pmb{\zeta}^{*}(\rho) = 0}C_{\rho}y^{1 - \rho/2} \qquad y \rightarrow 0, \label{a}
\eeq

\vspace{.27 cm}

 is $f$ of  rapid decay  at the cusp $y \rightarrow \infty$,
plus   some degree of  smoothness  of the function $f(x,y)$,  and   suitable  $y \rightarrow \infty$ growing conditions
for $\Delta f$, (where $\Delta : = y^{2}(\partial^{2}_{x} + \partial^{2}_{y})$ is the hyperbolic Laplacian).
 We make this precise, and derive a sufficient condition for (\ref{a}) to occur.

\vspace{.35 cm}

The starting point is the Rankin-Selberg integral

\vspace{.27 cm}

\beq
\pmb{I}(s) = \int_{\pmb{\mathcal{D}}} dx dy \, y^{-2} \pmb{E}_{s}(z) f(z). \label{I}
\eeq

\vspace{.27 cm}

Since the integral function $\pmb{I}(s)$ inherits analytic properties of $\pmb{E}_{s}(z)$, $\pmb{I}(s)$
 has a meromorphic continuation with   poles in $s = 1$, and $s = \rho/2$, with $\rho$'s such that $\pmb{\zeta}^{*}(\rho) = 0$.
Define $\Theta : = Sup\{\Re(\rho)| \pmb{\zeta}^{*}(\rho) = 0 \}$, $1/2 \le \Theta < 1$,
then $\pmb{I}(s) - \frac{C}{s - 1}$ is defined on $\Re(s) > \Theta/2$.

\vspace{.35 cm}

Since  $\pmb{I}(s) = \pmb{\mathcal{M}}[y^{-1}\pmb{a}_{0}(y)](s)$, a way  to obtain (\ref{a}) is to use an inverse Mellin transform argument [Ve].
The  Mellin inverse-transform  of $\pmb{I}(s)$  is

\vspace{.27 cm}

\beq
\pmb{\mathcal{M}}^{-1}[\pmb{I}(s)](y) = \frac{1}{2\pi i}\int_{\sigma - i\infty}^{\sigma + i\infty}ds \, y^{-s}\pmb{I}(s) = \frac{y^{- \sigma}}{2 \pi i} \int_{-\infty}^{\infty}dt \, y^{-it}\pmb{I}(\sigma + it), \label{M}
\eeq

\vspace{.27 cm}

wherever  $\pmb{I}(\sigma + i t )$ falls off as $o(1/t)$ for $t \rightarrow \pm \infty$.

If $f(z)$ is twice differentiable, then one can use $\Delta \pmb{E}_{s}(z) = s(1 - s)\pmb{E}_{s}(z)$
and by integration by parts one finds
\beq
\pmb{I}(s) =  \frac{1}{s(s - 1)}\int_{\pmb{\mathcal{D}}} dx dy \, y^{-2} \pmb{E}_{s}(z) \Delta f(z). \label{I2}
\eeq
This shows that $I(\sigma + it)$ falls off as $t^{-2}$ for $t \rightarrow \pm \infty$, whenever the integral
r.h.s. of (\ref{I2})  is convergent.
 For our purposes one has to check that this integral is convergent in $\sigma = \frac{\Theta}{2} + \epsilon$.
For $y \rightarrow \infty$, the Eisenstein series goes as $E_{z}(s) \sim y^{s} + \frac{\pmb{\zeta}^{*}(2s -1)}{\pmb{\zeta}^{*}(2s)}y^{1 - s}$,
and since  $1/4 <  \sigma = \frac{\Theta}{2} + \epsilon < 1/2$, indeed  $E_{\frac{\Theta}{2} + \epsilon}(z) \sim y^{1 - \frac{\Theta}{2} - \epsilon}$.
Thus the integral in (\ref{I2}) is convergent if $\Delta f(z)$ respects and upper bound for its polynomial growing $y \rightarrow \infty$,
 namely $\Delta f(z) \lesssim O(y^{1/4})$.

Alltogether, we  have the following sufficient condition:

\vspace{.3 cm}

\begin{propo} \label{p-rapid}
Given $f = f(x,y)$ a  modular invariant function  of rapid decay $y \rightarrow \infty$. If $f$ is twice differentiable and $\Delta f \lesssim O(y^{1/4})$
for $y \rightarrow \infty$, then the following holds true

\vspace{.27 cm}

\beq
\pmb{a}_{0}(y) \sim C + O(y^{1 - \frac{\Theta}{2}}) \qquad y \rightarrow 0,
\eeq

\vspace{.27 cm}

with $\Theta : = Sup\{\Re(\rho)| \pmb{\zeta}^{*}(\rho) = 0 \}$.
\end{propo}

\vspace{.3 cm}

\subsection{Rate of uniform distribution of long horocycles}

For the rate of uniform distribution of  horocycles  $\pmb{\mathcal{H}}_{y} : = (\mathbb{R} + i y) /\pmb{\Gamma}_{\infty} \subset \pmb{\mathcal{D}}$,
 in the modular surface $\pmb{\mathcal{D}} \simeq   \pmb{\Gamma}  \backslash \pmb{\mathcal{H}}$, one can prove  that

\vspace{.27 cm}

\beq
\frac{L(\pmb{\mathcal{H}}_{1/y} \cap \pmb{\mathcal{U}}) }{L(\pmb{\mathcal{H}}_{1/y})} \sim  \frac{A(\pmb{\mathcal{U}})}{A(\pmb{\mathcal{D}})} + O(y^{1/2}),   \qquad y \rightarrow 0 \label{h}
\eeq

\vspace{.27 cm}

for every open set $\pmb{\mathcal{U}} \subset \pmb{\mathcal{D}}$. $L$ indicates hyperbolic length, ($L(\gamma) = \int_{\gamma} y^{-1}\sqrt{dx^{2} + dy^{2}}$ for
a given curve $\gamma \subset \pmb{\mathcal{H}}$), and $A$ hyperbolic area $A(\pmb{\mathcal{U}}) = \int_{\pmb{\mathcal{U}}} dx dy y^{-2}$.
Eq. (\ref{h}) shows that for every open set $\pmb{\mathcal{U}}$ contained in $\pmb{\mathcal{D}}$,
the portion of horocycle $\pmb{\mathcal{H}}_{y}$ contained in $\pmb{\mathcal{U}}$ in the limit $y \rightarrow 0$
tends to become proportional to the ratio between the area $A(\pmb{\mathcal{U}})$ of $\pmb{\mathcal{U}}$,  and the area $A(\pmb{\mathcal{D}}) = \pi/3$ of $\pmb{\mathcal{D}}$.

\vspace{.35 cm}

The missing presence of $\Theta = Sup\{\Re(\rho)| \pmb{\zeta}^{*}(\rho) = 0 \}$  and thus the missing link with the Riemann hypothesis
 in the error estimate of (\ref{h})   is due to the fact that some of the arguments used
to prove proposition (\ref{p-rapid}) do not go through in the present case.
 In fact, one has

\vspace{.27 cm}

 \beq
 \frac{L(\pmb{\mathcal{H}}_{1/y} \cap \pmb{\mathcal{U}}) }{L(\pmb{\mathcal{H}}_{1/y})} = \int_{0}^{1} dx \, \, \pmb{\chi}_{\pmb{\mathcal{U}}}(x + iy),
 \eeq

\vspace{.27 cm}

where $\pmb{\chi}_{\pmb{\mathcal{U}}}(z)$ is the characteristic function of $\pmb{\mathcal{U}} \subset \pmb{\mathcal{D}}$.
 Also, by using the Rankin-Selberg method

\vspace{.27 cm}

 \beq
I_{\pmb{\chi}}(s) : = \pmb{\mathcal{M}}\left(\frac{1}{y}\frac{L(\pmb{\mathcal{H}}_{1/y} \cap \pmb{\mathcal{U}}) }{L(\pmb{\mathcal{H}}_{1/y})}\right)(s) = \int_{\pmb{\mathcal{D}}} dx dy \, y^{-2} \pmb{\chi}_{\pmb{\mathcal{U}}}(z) \pmb{E}_{z}(s).
\eeq

\vspace{.27 cm}

Since $\pmb{\chi}_{\pmb{\mathcal{U}}}(z)$ is not smooth, one cannot use   the Laplacian $\pmb{\Delta}$ argument,
 as it was done for deriving  proposition \ref{p-rapid}.
Thus, the  inverse-Mellin  argument does not go through, and there is no
connection between the rate of uniform distribution of long horocycles in the
modular surface $ \pmb{\Gamma} \backslash \pmb{\mathcal{H}}$  and the Riemann hypothesis.

\vspace{.3 cm}

\section{\Large{\textbf{Modular functions of polynomial growth}}}\label{SPpolynomial}

\vspace{.3 cm}

For a modular invariant function $f$ of polynomial growth at the cusp

$$f(z) \sim \varphi(y)  + o(y^{-N}), \qquad y \rightarrow \infty  \qquad \forall N > 0 $$

where
\beq
 \varphi(y) : = \sum_{i=1}^{l}\frac{c_i}{n_{i}!} \, y^{\alpha_i}\, \log^{n_i}y, \label{varphi}
\eeq
and

$$\alpha_i, c_i \in \mathbb{C},  \Re(\alpha_{i}) < 1/2, n_i \in \mathbb{N}_{\ge 0}.$$

\vspace{.3 cm}
 Zagier  [Za2] has obtained  analytic continuation and functional equation of the following Rankin-Selberg
 integral transform
\eqn
\pmb{R}^{*}(f,s) :&=& \pmb{\zeta}^{*}(2s) \int_{0}^{\infty} dy y^{s - 2}(\pmb{a}_{0}(y) - \varphi(y) ) \nn \\
&=& \sum_{i=1}^l c_i \left( \frac{\pmb{\zeta}^{*}(2s)}{(1 - s - \alpha_i)^{n_i + 1}}
 + \frac{\pmb{\zeta}^{*}(2s - 1)}{(s - \alpha_i )^{n_i + 1}}  + \frac{{\rm{ entire \,\, function \,\, of}} \,\, s}{s(s - 1)}    \right). \nn \\ \label{eq-RSZpolar2}
\feqn

\vspace{.35 cm}

Eq. (\ref{eq-RSZpolar2}) is obtained  in [Za2] by a method which in terms of the following diagram

 \vspace{.3 cm}

 \beq
\begin{array}[c]{ccc}
\langle \pmb{\Theta}_{t}(z), f(z) \rangle_{ \pmb{\Gamma} \backslash \pmb{\mathcal{H}}}    &      \overset{\pmb{\mathcal{M}}_{t} }{\longrightarrow }& \langle \pmb{E}^{*}_{s}(z), f(z) \rangle_{ \pmb{\Gamma} \backslash \pmb{\mathcal{H}}}   \\
& & \\
\downarrow Unfolding  &                  &   \downarrow   Unfolding  \\
&  & \\
\langle \pmb{\vartheta}_{t}(y), \pmb{a}_{0}(y) \rangle_{U \backslash \pmb{\mathcal{H}} } &   \overset{\pmb{\mathcal{M}}_{t} }{\longrightarrow } &  \pmb{\zeta}^{*}(2s)\langle y^{s}, \pmb{a}_{0}(y) \rangle_{U \backslash \pmb{\mathcal{H}} },
\end{array} \label{magic}
\eeq

 \vspace{.35 cm}

 corresponds in considering the Rankin-Selberg integral in the right upper vertex of  diagram \ref{magic},
   albeit with a regularization in
 the integration in the presence of a  cutoff $T > 1$,   $\pmb{\mathcal{D}}_T = \{z \in \pmb{\mathcal{D}}| y < T, \}$.
 This truncation allows to apply a  version of the unfolding trick devised for truncated domains $\pmb{\mathcal{D}}_T$, and to
 move along the  right column of this diagram. The obtained unfolded
 $T$-dependent quantity comprises several terms, and  a  careful analysis of the $T \rightarrow \infty$ limit [Za2]
 allows to extract information on $\pmb{R}^{*}(f,s) = \pmb{\zeta}^{*}(2s)\langle y^{s}, \pmb{a}_{0}(y) \rangle_{U \backslash \pmb{\mathcal{H}} }$,
 in the lower right corner of the diagram \ref{magic}.
 This leads  to prove  equation (\ref{eq-RSZpolar2}) for the meromorphic continuation of $\pmb{R}^{*}(f,s)$,
 plus additional results on functional equation for the Rankin-Selberg transform $\pmb{R}^{*}(f,s)$ [Za2].

 \vspace{.35 cm}

Here we employ an alternative method which leads us to prove   (\ref{eq-RSZpolar2}). This method allows us
to obtain results on the long horocycle average of functions with growing conditions
given in (\ref{varphi}), i.e. functions in $\pmb{\mathcal{C}}_{TypeII}$.
Our method comprises the following  two steps in the diagram

\vspace{.5 cm}

\beq
\begin{array}[c]{ccc}
\langle \pmb{\Theta}_{t}(z), f(z) \rangle_{ \pmb{\Gamma} \backslash \pmb{\mathcal{H}}}    &    &    \\
& & \\
\downarrow Unfolding  &                     &  \\
&  & \\
\langle \pmb{\vartheta}_{t}(y), \pmb{a}_{0}(y) \rangle_{U \backslash \pmb{\mathcal{H}} } &   \overset{\pmb{\mathcal{M}}_{t} }{\longrightarrow } &  \pmb{\zeta}^{*}(2s)\langle y^{s}, \pmb{a}_{0}(y) \rangle_{U \backslash \pmb{\mathcal{H}} }  = \pmb{R}^{*}(f,s).
\end{array} \label{magic2bis}
\eeq

\vspace{.5 cm}

The advantage of this route is that it does not require regularization (truncations)  of the domain of integration $\pmb{\mathcal{D}}$.
In order to perform the unfolding of the integration domain in the modular integral

\vspace{.27 cm}

\beq
\langle \pmb{\Theta}_{t}(z), f(z) \rangle_{ \pmb{\Gamma} \backslash \pmb{\mathcal{H}}} = \int_{\pmb{\mathcal{D}}} dx dy \, y^{-2} f(x,y) \, \sum_{(m,n) \in \mathbb{Z}^2 \setminus \{ 0 \}} e^{- \frac{\pi}{y}|m z +  n|^{2}}, \label{TheraInt}
\eeq

\vspace{.27 cm}

 from the decomposition for the  theta series $\pmb{\Theta}_{t}(z)$

\vspace{.27 cm}

 \beq
 \pmb{\Theta}_{t}(z) = \sum_{(m,n) \in \mathbb{Z}^2 \setminus \{ 0 \}}e^{-\pi t \frac{(mx + n)^{2} + m^2 y^2 }{y}} =
  1 + \pmb{\vartheta}(t/y) + \sum_{\gamma \in \pmb{\Gamma}_{\infty} \backslash  \pmb{\Gamma}^{'}} \pmb{\vartheta}_{t}(1/\Im(\gamma(z))), \nn
 \eeq

\vspace{.27 cm}

 one uses contributions from the third term on the r.h.s.,
 where $\pmb{\Gamma}^{'} : =  \pmb{\Gamma} \, \, \setminus \, \, \{ \mathbb{I} \}$ is the set of modular transformations minus the identity $\mathbb{I}$.
 Each term in this  series has the form

\vspace{.27 cm}

\beq
 \pmb{\vartheta}_{t}(1/\Im(\gamma(z))) = \sum_{r \ne 0}e^{-\pi \frac{r^2}{y} ((c x + d)^{2} + c^2 y^2)} \qquad c, d \in \mathbb{Z}, c \ne 0, (c,d) = 1 \nn
 \eeq

\vspace{.27 cm}

and  corresponds  to the $m \ne 0$ subseries in (\ref{TheraInt}), whose terms decay exponentially for $y \rightarrow \infty$.
 Thus by dominate convergence theorem one can
unfold the modular integral $\langle \pmb{\Theta}_{t}(z), f(z) \rangle_{ \pmb{\Gamma} \backslash \pmb{\mathcal{H}}} $
in the left upper entry of (\ref{magic2bis})
and prove the vertical arrow connecting the left upper entry with  the left lower entry $\langle \pmb{\vartheta}_{t}(y), \pmb{a}_{0}(y) \rangle_{U \backslash \pmb{\mathcal{H}} }$.

\vspace{.35 cm}

The unfolding  trick is doable  without using a truncated domain,
since the integral in the left upper corner of the diagram  is convergent,  under the assumptions
$\Re(\alpha_{i}) < 1/2$, for the growing term $\varphi(y)$ in (\ref{varphi}). Indeed, by Poisson summation
one can check that $\pmb{\Theta}_{t}(z) \sim  \sqrt{y}$ for $y \rightarrow \infty$.
 Moreover, $\pmb{\Theta}_{t}(z)$  has  series representation convergent  for every $t > 0$.
 Thus we have the following proposition for Theta-unfolding of a modular invariant  function $f$
with growing conditions in $\pmb{\mathcal{C}}_{TypeII}$ (\ref{poly}):

\vspace{.3 cm}

\begin{propo} \label{ThetaUnfPoly}

Let $f = f(x,y)$ a modular invariant function of polynomial growth at the cusp $y \rightarrow \infty$

\beq
  f(x,y) \sim    \sum_{i=1}^{l}\frac{c_i}{n_{i}!} \, y^{\alpha_i}\, \log^{n_i}y +  o(y^{-N}), \qquad \forall N \ge 0, \qquad y \rightarrow \infty \nn
\eeq
with
$$\alpha_i, c_i \in \mathbb{C},  \Re(\alpha_{i}) < 1/2, n_i \in \mathbb{N}_{\ge 0}.$$

Then, the following Theta-unfolding relation holds true

\beq
\int_{\pmb{\mathcal{D}}} dx dy \, y^{-2} f(x,y) \pmb{\Theta}_{t}(z) = \int_{0}^{\infty}dy \, y^{-2} \pmb{a}_{0}(y) \pmb{\vartheta}_{t}(y).
\eeq
\end{propo}

\vspace{.3 cm}

Proposition \ref{ThetaUnfPoly} states that the  vertical arrow in diagram (\ref{magic2bis})
holds true for modular functions of polynomial growth class $\pmb{\mathcal{C}}_{TypeII}$.

\vspace{.35 cm}

The horizontal arrow in diagram (\ref{magic2bis}) indicates that  due to the relation between the
functions

\vspace{.27 cm}

\beq
\pmb{i}(t) : = \langle \pmb{\vartheta}_{t}(y), \pmb{a}_{0}(y) \rangle_{U \backslash \pmb{\mathcal{H}}} \label{i}
\eeq

\vspace{.27 cm}

and the function $\pmb{\zeta}^{*}(2s)\langle y^{s}, \pmb{a}_{0}(y) \rangle_{U \backslash \pmb{\mathcal{H}} }$
through Mellin transform, knowledge of the $t \rightarrow \infty$ and $t \rightarrow 0$ asymptotics
for $\pmb{i}(t)$ implies knowledge of the meromorphic continuation with orders and locations of  poles
of the function  $\pmb{\zeta}^{*}(2s)\langle y^{s}, \pmb{a}_{0}(y) \rangle_{U \backslash \pmb{\mathcal{H}} }$
of complex variable $s$.
We therefore prove $\pmb{i}(t)$ asymptotics in two steps, by the two following lemmas.


\vspace{.35 cm}

\begin{lemma2} {\rm{\textbf{\ref{lemma1}.}}} \label{lemma1b} Let $f =  f(x,y)$ a modular invariant functions with growing conditions
 as  in proposition \ref{ThetaUnfPoly}.  Let  $C_0 : = < 1, f >_{ \pmb{\Gamma} \backslash \pmb{\mathcal{H}}}$,
  its integral over the  fundamental domain $\pmb{\mathcal{D}}$,
 and let $\pmb{a}_{0}(y)$ be  the $f$ constant term.

 Then,  the following relation holds true

\vspace{.27 cm}

\beq
\langle \pmb{\vartheta}_{t}(y), \pmb{a}_{0}(y) \rangle_{U \backslash \pmb{\mathcal{H}}} = \frac{1}{t}\langle \pmb{\vartheta}_{1/t}(y), \pmb{a}_{0}(y) \rangle_{U \backslash \pmb{\mathcal{H}}} + \frac{C_{0}}{t} - C_{0} \nn
\eeq

\vspace{.27 cm}

\end{lemma2}

\begin{proof}
By double Poisson summation one finds  $\pmb{\Theta}_{t}(z) = \frac{1}{t}\pmb{\Theta}_{1/t}(z) - \frac{1}{t} - 1$.
 The thesis then  follows by  applying the Theta-unfolding in proposition \ref{ThetaUnfPoly},
 which corresponds of using the left column in diagram  \ref{magic2bis}.

\end{proof}


\vspace{.3 cm}


By previous lemma, we are now in the position of proving the following lemma on
the asymptotics $t \rightarrow \infty$ and $t \rightarrow 0$ of the function $\pmb{i}(t)$ defined
by (\ref{i}), (which appears in the left lower entry of  diagram \ref{magic2bis}):

\vspace{.35 cm}

\begin{lemma2} {\rm{\textbf{\ref{propolemma2}.}}}     \label{lemma2}
Let $f =  f(x,y)$ a modular invariant function  with polynomial behavior at the cusp

\vspace{.27 cm}

\beq
f(z) \sim \sum_{i=1}^{l}\frac{c_i}{n_{i}!} \, y^{\alpha_i}\, \log^{n_i}y  + o(y^{-N}), \qquad y \rightarrow \infty  \qquad \forall N > 0  \nn
\eeq

\vspace{.27 cm}

where $\alpha_i, c_i \in \mathbb{C}$,  $\Re(\alpha_{i}) < 1/2$, $n_i \in \mathbb{N}_{\ge 0}$.

\vspace{.2 cm}

Then, for the function $\pmb{i}(t) : = \langle \pmb{\vartheta}_{t}(y), \pmb{a}_{0}(y) \rangle_{U \backslash \pmb{\mathcal{H}}}$ the following
asymptotics hold true
\eqn
&i)& \qquad \pmb{i}(t) \sim \sum_{i=1}^{l}\frac{c_i}{n_{i}!}  \pmb{\zeta}^{*}(2\alpha_{i} - 1 ) t^{\alpha_{i} - 1}\, \log^{n_i} t + O\left(t^{A  - 1}\log^{N - 1} t   \right),   \qquad  t \rightarrow \infty   \nn \\
&ii)& \qquad  \pmb{i}(t) \sim  - C_{0} +  \frac{C_{0}}{t} - \sum_{i=1}^{l}\frac{c_i}{n_{i}!}    \pmb{\zeta}^{*}(2\alpha_{i} - 1) t^{- \alpha_{i}} \log^{n_i} t
 + O\left(t^{-a}\log^{n - 1} t   \right)  ,  \qquad   t \rightarrow 0 \nn
\feqn
where   $A : = max\{ \Re(\alpha_i) \}$, $a : = min\{ \Re(\alpha_i) \}$, $N : =   max\{ n_i \}$,  $n : =  min\{ n_i \}$, and
$$C_{0} = \int_{\pmb{\mathcal{D}}}dx dy y^{-2} f(z).   $$
\end{lemma2}

\vspace{.3 cm}

\begin{proof}
We start by proving i),   $\pmb{i}(t)$ is the following integral function
\beq
\pmb{i}(t) = \int_{0}^{\infty}dy y^{-2} \pmb{a}_{0}(y)\sum_{r \in \mathbb{Z} \backslash \{0 \}}e^{-\pi r^{2}\frac{t}{y}}, \nn
\eeq
by change of integration variable $y \rightarrow t  y$ one finds

\beq
\pmb{i}(t) = \frac{1}{t} \int_{0}^{\infty}dy y^{-2} \pmb{a}_{0}(y t)\sum_{r \in \mathbb{Z} \backslash \{0 \}}e^{-\pi \frac{r^{2}}{y}}. \nn
\eeq
Therefore  for $t \rightarrow \infty$

\eqn
\pmb{i}(t) &\sim& \sum_{i=1}^{l}\frac{c_i}{n_{i}!}   t^{\alpha_i - 1}  \int_{0}^{\infty}dy y^{-2 + \alpha_{i}}\left( \log y + \log t   \right)^{n_{i}}  \sum_{r \in \mathbb{Z} \backslash \{0 \}}e^{-\pi \frac{r^{2}}{y}}, \nn \\
&\sim&  \sum_{i=1}^{l}\frac{c_i}{n_{i}!}  \pmb{\zeta}^{*}(2\alpha_{i} - 1 ) t^{\alpha_{i} - 1}\, \log^{n_i} t + O\left(t^{A  - 1} \log^{N - 1}y   \right). \nn
\feqn

In order to prove  ii),   we use lemma \ref{lemma1b} which allows  to rewrite $\pmb{i}(t)$ in the following form
\eqn
\pmb{i}(t) &=& \frac{1}{t}\int_{0}^{\infty}dy y^{-2} \pmb{a}_{0}(y)\pmb{\vartheta}_{1/t}(y) + \frac{C_{0}}{t} - C_{0} \nn \\
   &=&  \int_{0}^{\infty}dy y^{-2} \pmb{a}_{0}(y/t) \sum_{r \in \mathbb{Z} \backslash \{0 \}} e^{-\pi \frac{r^{2}}{y}} + \frac{C_{0}}{t} - C_{0}, \nn
\feqn
also for $t \rightarrow 0$
\eqn
&& \int_{0}^{\infty}dy y^{-2} \pmb{a}_{0}(y/t) \sum_{r \in \mathbb{Z} \backslash \{0 \}} e^{-\pi \frac{r^{2}}{y}}  \nn \\
&\sim&
\sum_{i=1}^{l}\frac{c_i}{n_{i}!}   \int_{0}^{\infty}dy y^{-2 - \alpha_i} \left( \log y - \log t   \right)^{n_{i}}    \sum_{r \in \mathbb{Z} \backslash \{0 \}} e^{-\pi \frac{r^{2}}{y}} \nn \\
&\sim&
- \sum_{i=1}^{l}\frac{c_i}{n_{i}!}    \pmb{\zeta}^{*}(2\alpha_{i} - 1) t^{- \alpha_{i}} \log^{n_i} t
 + O\left(t^{-a}\log^{n - 1} t   \right)  \nn
\feqn
\end{proof}

\vspace{.3 cm}

In order to prove   Zagier  theorem \ref{Zpol}
 on the analytic continuation of the Rankin-Selberg transform,  from  lemma \ref{lemma2} and  from the lower row of diagram \ref{magic},
we also need  proposition \ref{propomellin1} on standard properties of Mellin transforms.
Due to the lower row in diagram \ref{magic2},
by applying  proposition  \ref{propomellin1} on the asymptotics  in lemma \ref{lemma2}, we  obtain analytic continuation of
the Rankin-Selberg transform, as in theorem
\ref{Zpol}.

\vspace{.3 cm}

From lower row of diagram \ref{magic2}, lemma \ref{lemma2} and proposition \ref{propomellin1},
we also prove  the following theorem on the asymptotic of the long horocycle average
of a modular function in $\pmb{\mathcal{C}}_{TypeII}$:

\vspace{.3 cm}

\begin{theorem3}{\rm{\textbf{\ref{thpoly2}.}}}       \label{thpoly4}
Let  $f = f(x,y)$ a modular invariant function   of  polynomial growth at the cusp

 $$f(x,y) \sim \sum_{i=1}^{l}\frac{c_i}{n_{i}!} \, y^{\alpha_i}\, \log^{n_i}y  + o(y^{-N}), \qquad y \rightarrow \infty  \qquad \forall N > 0 $$

where $\alpha_i, c_i \in \mathbb{C}$,  $\Re(\alpha_{i}) < 1/2$, $n_i \in \mathbb{N}_{\ge 0}$.

\vspace{.3 cm}

The long length limit of the $f$ horocycle average has the following asymptotic

\beq
 \pmb{a}_{0}(y) \sim C_{0} + \sum_{\pmb{\zeta}^{*}(\rho) = 0} C_{\rho} y^{1 - \frac{\rho}{2}} +
  \sum_{i=1}^{l}\frac{c_i}{n_{i}!}    \frac{\pmb{\zeta}^{*}(2\alpha_i - 1)}{\pmb{\zeta}^{*}(2\alpha_i )}y^{1 - \alpha_i }\log^{n_i} y + O\left(y^{1 - A}\log^{n - 1} y \right)     \qquad y \rightarrow 0, \nn
 \eeq
where   $A : = max\{ \Re(\alpha_i) \}$,  $n : =  min\{ n_i \}$, and
$$C_{0} = \frac{3}{\pi}\int_{\pmb{\mathcal{D}}}dx dy \, y^{-2} f(z).   $$
\end{theorem3}

\vspace{.3 cm}


\section{\Large{\textbf{Modular functions of exponential growth}}}\label{SPexponential}

\vspace{.3 cm}

We now turn to discuss modular invariant functions in the class of growing conditions
$\pmb{\mathcal{C}}_{Heterotic}$, defined in (\ref{cstring}).
Proofs are obtained by  using same methods we employed in previous sections
for the  $\pmb{\mathcal{C}}_{Type II}$ case, i.e. by following  the arrows in diagram \ref{magic2bis}.

\vspace{.35 cm}

We start by proving a bound on the growing of the long horocycle average for
a modular function in $\pmb{\mathcal{C}}_{Heterotic}$. Some  of the ideas contained  in the proof of theorem \ref{thbound}
are taken  from [KS].

\vspace{.3 cm}

\begin{theorem3}\label{thExpBound} {\rm{\textbf{\ref{thbound}.}}}
Let $f = f(x,y)$ be a modular invariant function with the following growing condition

\vspace{.27 cm}
\beq
  f(x,y) \sim y^{\alpha} e^{2\pi i \kappa x}e^{\pi \beta y} \qquad y \rightarrow \infty \qquad  \kappa \in \mathbb{Z}  \, \, \setminus \, \, \{ 0 \}, \qquad \beta < 1, \qquad \alpha \in \mathbb{C}, \, \,  \Re(\alpha) < 1/2. \nn
\eeq

\vspace{.27 cm}

 Then the  growth of  the long horocycle average $\pmb{a}_{0}(y)$  satisfies
 the following bound
 $$  \pmb{a}_{0}(y) \lesssim o(e^{C/y}) \qquad  y \rightarrow 0, \qquad \forall C \in \mathbb{C}, \, \,  \Re(C) > 0.$$
\end{theorem3}

\vspace{.35 cm}

\begin{proof}

We consider the  following  Theta-integral on the modular domain $\pmb{\mathcal{D}}$

\vspace{.27 cm}

\beq
\pmb{I}(t) : = \int_{\pmb{\mathcal{D}}} dx dy y^{-2} f(z) \sum_{(m,n) \in \mathbb{Z}^{2} \backslash \{ 0 \}}e^{- \frac{\pi t}{y}|m x + n|^{2}}, \label{Itheta}
\eeq

\vspace{.27 cm}

which corresponds to Petersson inner product of the theta series $\pmb{\Theta}_{t}(z)$ with the function $f$
which appears  in the left upper entry of diagram \ref{magic2bis}.
Due to the $f$ growing conditions for $y \rightarrow \infty$, the function $\pmb{I}(t)$, for small $t$, has to be understood
as the result of an integration over $\pmb{\mathcal{D}}$, with integration along the real axis performed first.
In fact, the modular integral is only conditionally convergent for $z \rightarrow i \infty$.

\vspace{.35 cm}

We employ the following decomposition for the theta series $\pmb{\Theta}_{t}(z)$

\vspace{.27 cm}

\beq
\pmb{\Theta}_{t}(z)  =  \sum_{(m,n) \in \mathbb{Z}^{2} \backslash \{ 0 \}}e^{- \frac{\pi t}{y}|m z + n|^{2}} = \sum_{\mathbb{Z} \backslash \{ 0 \}}e^{- \pi t \frac{r^{2}}{y}} +
\sum_{\pmb{\Gamma}_{\infty}  \backslash \pmb{\Gamma}^{'}}\sum_{\mathbb{Z} \backslash \{ 0 \}}  e^{- \pi t \frac{r^{2}}{\Im(\gamma(z))}}, \nn
\eeq

\vspace{.27 cm}

where $\pmb{\Gamma}^{'} = \pmb{\Gamma} \, \backslash \, \{ \mathbb{I} \}$ is  the modular group $\pmb{\Gamma}$ minus the identity $\mathbb{I}$.

One has for $\gamma \in \pmb{\Gamma}_{\infty}  \backslash \pmb{\Gamma}^{'}$

\vspace{.27 cm}

\beq
e^{- \pi t \frac{r^{2}}{\Im(\gamma(z))}} \sim e^{-  \pi t  m^{2} y },   \qquad y \rightarrow \infty, \nn
\eeq

\vspace{.27 cm}

with $m = c r$ and $c \ne 0$ is  the third entry of the modular transformation $\gamma$.
Modular transformations in $\gamma \in \pmb{\Gamma}_{\infty}  \backslash \pmb{\Gamma}^{'}$ allow to unfold integration
domain $\pmb{\mathcal{D}} \simeq  \pmb{\Gamma} \backslash \pmb{\mathcal{H}}$ in $\pmb{I}(t)$ into the half-infinite strip $\pmb{\Gamma}_{\infty} \backslash \pmb{\mathcal{H}}$.
From Lebesgue dominated convergence theorem, for $t > 1$ one finds

\vspace{.27 cm}

\beq
\int_{\pmb{\mathcal{D}}} dx dy y^{-2} f(z) \sum_{\pmb{\Gamma}_{\infty}    \backslash \pmb{\Gamma}^{'}}\sum_{\mathbb{Z} \backslash \{ 0 \}}  e^{- \pi t \frac{r^{2}}{\Im(\gamma(z))}}
=  \sum_{\pmb{\Gamma}_{\infty}    \backslash \pmb{\Gamma}^{'}}\sum_{\mathbb{Z} \backslash \{ 0 \}}\int_{\pmb{\mathcal{D}}} dx dy y^{-2} f(z)e^{- \pi t \frac{r^{2}}{\Im(\gamma(z))}}.
\eeq

\vspace{.27 cm}

This leads to the following  Theta-unfolding relation ($t > 1$)

\vspace{.27 cm}

\beq
\sum_{\pmb{\Gamma}_{\infty}    \backslash \pmb{\Gamma}^{'}}\sum_{\mathbb{Z} \backslash \{ 0 \}}\int_{\pmb{\mathcal{D}}} dx dy y^{-2} f(z)e^{- \pi t \frac{r^{2}}{\Im(\gamma(z))}} = \int_{0}^{\infty} dy y^{-2} \int_{-1/2}^{1/2} dx  f(x,y) \sum_{\mathbb{Z} \backslash \{ 0 \}}e^{- \pi t \frac{r^{2}}{y}}. \label{ThetaUnfExp}
\eeq

\vspace{.27 cm}

Therefore, for $t > 1$, the following unfolding relation holds

\vspace{.27 cm}

\beq
\pmb{I}(t) = \int_{0}^{\infty} dy y^{-2} \pmb{a}_{0}(y) \sum_{\mathbb{Z} \backslash \{ 0 \}}e^{- \pi t \frac{r^{2}}{y}}. \label{It}
\eeq

\vspace{.27 cm}

Moreover, one can prove that the function $\pmb{I}(t)$ in her original incarnation (\ref{Itheta}),   is analytic on a strip  $t \in (0, \infty) \times (- \delta_{\beta}, \delta_{\beta}) \subset
\mathbb{C}$, where $\delta_{\beta} : =  1 - \beta > 0$.
Proof of this statement   follows by  Poisson summation

\vspace{.27 cm}

\beq
\pmb{I}(t) =  \frac{1}{\sqrt{t}}\int_{\pmb{\mathcal{D}}} dx dy y^{-3/2} f(z) \sum_{(m,n) \in \mathbb{Z}^{2} \backslash \{ 0 \}} e^{- \pi y \left(\frac{m^{2}}{t} + n^{2}t   \right)}e^{2\pi i m n x},
\eeq

\vspace{.27 cm}

and by the growing assumption we make on $f(x,y)$ for  $y \rightarrow \infty$.

\vspace{.3 cm}

Due to analyticity of the l.h.s. in (\ref{ThetaUnfExp}) on the strip $t \in (0, \infty) \times (- \delta_{\beta}, \delta_{\beta}) \subset
\mathbb{C}$, where $\delta_{\beta}  =  1 - \beta > 0$, the r.h.s. cannot be divergent on this strip.
This rules out the following behavior

\beq
\pmb{a}_{0}(y) = \int_{0}^{1}dx f(x,y) \sim e^{C/y} \qquad y \rightarrow 0 \qquad C \in \mathbb{C}, \, \, \Re(C) > 0, \nn
\eeq
since such a growing condition would make  the integral function in the r.h.s. of (\ref{ThetaUnfExp}) to diverge
for $0 <  t < \Re(C)$.
\end{proof}

\vspace{.3 cm}

As  already  remarked, string theory suggests a much stronger result than theorem \ref{thExpBound},
 namely  that in the $y \rightarrow 0$ limit  $\pmb{a}_{0}(y)$
be convergent and to have  asymptotic
 as in theorem \ref{thpoly4}.

This leads to the following open question:

\vspace{.27 cm}

\begin{open2}  {\rm{\textbf{\ref{Qexp}.}}}
 Given $f = f(x,y)$ modular invariant function in the class $\pmb{\mathcal{C}}_{Heterotic}$ (\ref{cstring}),
 prove or disprove that  the following asymptotic holds true
$$ \pmb{a}_{0}(y) \sim C_{0} + \sum_{\pmb{\zeta}^{*}(\rho) = 0} C_{\rho} y^{1 - \frac{\rho}{2}} + \frac{\pmb{\zeta}^{*}(2\alpha - 1)}{\pmb{\zeta}^{*}(2\alpha)}y^{1 - \alpha} + o(y^{N}), \qquad  \forall N > 0, \qquad y \rightarrow 0, $$

$$C_{0} = \frac{3}{\pi}\int_{\pmb{\mathcal{D}}} dy \int dx y^{-2} f(z),   $$
where this integral is meant in the conditional sense, with integration along the real axis first  performed.
\end{open2}

\vspace{.27 cm}

Although we are not able to address the above question,
 we would like to end by adding few  remarks, which
 may be relevant to address the open problem \ref{Qexp}.
We consider the possibility that there may be some kind of rigidity in the
way the horocycle average can grow in the long length limit,
 for a  modular invariant function  $f$ with growing conditions
in $\pmb{\mathcal{C}}_{Heterotic}$.
Rigidity on the way $\pmb{a}_{0}(y)$ grows in the $y \rightarrow 0$ limit
under growing conditions on $f$  in  $\pmb{\mathcal{C}}_{Heterotic}$,
 may be contained in    proposition \ref{propmod}  below.

\vspace{.35 cm}

By using the following standard formulae for  transformations of the real and imaginary
part of $z \in  \pmb{\mathcal{H}}$ under a $SL_{2}(\pmb{Z})$ modular transformation $$\gamma(z) = \frac{az + b}{cz + d},$$ with $c \ne 0$,
one has

$$\Re(\gamma(z)) = \frac{a}{c}\frac{(x + b/a)(x + d/c) + y^2 }{(x + d/c)^{2} + y^2 },$$
$$\Im(\gamma(z)) = \frac{1}{c}\frac{y}{(x + d/c)^2 + y^2 }.$$

\vspace{.3 cm}

One can then prove the following proposition for modular functions in $\pmb{\mathcal{C}}_{Heterotic}$

\vspace{.3 cm}

\begin{propo2} {\rm{\textbf{\ref{propmod}.}}}
Given a $SL_{2}(\pmb{Z})$ invariant function $f$ which grows as $f(x,y) \sim e^{2\pi \beta y}e^{2\pi i \kappa x}$ for $y \rightarrow \infty$,
  $\kappa \in \mathbb{Z} \,  \backslash \,  \{0 \}$, Then

\vspace{.27 cm}

\beq
\pmb{a}_{0}(y) +  \sum_{r \in \mathbb{Z} \, \backslash \, \{ 0\}}\pmb{a}_{r}(y)e^{2\pi i r \frac{a}{c}} \sim e^{- 2 \pi i \kappa \frac{d}{c}} e^{2\pi \beta \frac{c^{2}}{y}},  \qquad y \rightarrow 0
\eeq

\vspace{.27 cm}

for every pair of Farey fractions $\frac{a}{c}$, $\frac{d}{c}$,  $a,c,d \in \mathbb{Z}$, $(a,c) = 1$, $|a| < c$, $(d,c) = 1$
$|d| < c$, $c > 0$.
 $\pmb{a}_{r}(y)$ are the Fourier modes in the expansion $f(x,y) = \sum_{r \in \mathbb{Z}} \pmb{a}_{r}(y)e^{2\pi i r x}$.

\end{propo2}

\vspace{.35 cm}

Perhaps   proposition \ref{propmod} together with theorem \ref{thbound}
turn out to  be sufficient    to address the  \textbf{open problem} \ref{Qexp}.
Another possibility is that question  \ref{Qexp} holds true in the following way.
It may be that all the
modular functions in the growing class  $\pmb{\mathcal{C}}_{Heterotic}$  split as the sum of a modular invariant  function
in $\pmb{\mathcal{C}}_{TypeII}$ plus a cusp function in $\pmb{\mathcal{C}}_{Heterotic}$, (a function whose constant term
$\pmb{a}_{0}(y)$ is identically  vanishing) \cite{Za3}.
We are not able to provide answers   on this latter possibility, which may give  an alternative way
that reconciles string theory suggestions with results in the automorphic function domain.

Finally, it could be as well possible that question raised in \ref{Qexp} has an answer
in the negative. This latter possibility would open some  interesting questions
in string theory, related to the emerging of a lack of  symmetry in the ultraviolet   between Type II  and Heterotic closed strings
asymptotics involving    very massive closed strings.

\section*{{\bf Acknowledgments}}
The Author thanks
 David Kazhdan, Gerard van der Geer and Don Zagier for enlighting  discussions.
The Author  thanks Anne-Marie Aubert for  reading the manuscript and
 providing useful suggestions and corrections, and
 Sergio Cacciatori for several  discussions, suggestions, and
collaboration on related topics. Thanks to
 Carlo Angelantonj, Shmuel Elitzur and Eliezer Rabinovici for
collaboration in string theory on  subjects  with   intersections
with the topic  of  this paper.
The Author  thanks
the Racah Institute of Physics at the Hebrew University of Jerusalem, the Theoretical Physics  Section
at the University of Milan Bicocca,
the ESI Schroedinger Center for Mathematical Physics in Vienna and the Theory Unit at
CERN for hospitality and support during  various  stages of this work. The Author  was supported
by a ``Angelo Della Riccia" visiting  fellowship at the University of Amsterdam.


\begin{thebibliography}{abc}





\bibitem[ACER]{ACER}   C.~Angelantonj, M.~Cardella, S.~Elitzur and E.~Rabinovici,
  ``Vacuum stability, string density of states and the Riemann zeta function,'' 
  JHEP {\bf 02(2011)}024, arXiv:1012.5091 [hep-th]

\bibitem[CC1]{CC}   S. Cacciatori,  M.~Cardella,
  ``Equidistribution rates,  closed string amplitudes, and the Riemann
  hypothesis,''
  JHEP {\bf 12 (2010)} 025
  arXiv:1007.3717 [hep-th].


\bibitem[CC2]{CC2}  S.~L.~Cacciatori and M.~A.~Cardella,  
  ``Uniformization, Unipotent Flows and the Riemann Hypothesis,''
  arXiv:1102.1201 [math.NT].
  
  \bibitem[CC3]{CC3}
  S.~L.~Cacciatori and M.~A.~Cardella,
  ``Eluding SUSY at every genus on stable closed string vacua,''
  arXiv:1102.5276 [hep-th].
  

\bibitem[C1]{C1}  M.~Cardella,
  ``A novel method for computing torus amplitudes for $\mathbb{Z}_{N}$
  orbifolds without the unfolding technique,'' 
  JHEP {\bf 0905} (2009) 010
  [arXiv:0812.1549 [hep-th]].


\bibitem[DS]{Dani} S. G. Dani and J. Smillie, Uniform distribution of horocycle orbits for Fuchsian groups. Duke Math. J. 51 (1984), 185194.

\bibitem[Fu]{F} H. Furstenberg, The Unique Ergodicity of the Horocycle Flow, Recent Advances in Topological Dynamics, A. Beck (ed.), Springer Verlag Lecture Notes, 318 (1972), 95-115.

\bibitem[He]{He}
 Gustav A. Hedlund,
``Fuchsian groups and transitive horocycles''
Duke Math. J. Volume 2, Number 3 (1936), 530-542.



\bibitem[KS]{KS}
  D.~Kutasov and N.~Seiberg,
``Number Of Degrees Of Freedom, Density Of States And Tachyons In String
  Theory And Cft,''
  Nucl.\ Phys.\  B {\bf 358}, 600 (1991).



\bibitem[R-S]{Rankin} R. Rankin, Contributions to the theory of Ramanujan's function $\tau(n)$ and
symilar arithmetical functions, Proc. Cambridge Philos. Soc. \textbf{35} (1939), 351-372. \\
 A. Selberg, Bemerkungen \`uber eine Dirichletsche reihe, die mit der theorie der modulformen nahe verbunden ist, Arch. Math.
Naturvid. \textbf{43} (1940), 47-50.

\bibitem[Ra]{Ratner1} M. Ratner,   Distribution rigidity for unipotent actions on homogeneous spaces.  Bull. Amer. Math. Soc. (N.S.) Volume 24, Number 2 (1991), 321-325. \\ M. Ratner,  Raghunathan's topological conjecture and distributions of unipotent flows.  Duke Math. J. Volume 63, Number 1 (1991), 235-280.



\bibitem[Sa]{Sa} P. Sarnak,
"Asymptotic behavior of periodic orbits of the horocycle flow and Eisenstein series." Comm. Pure Appl. Math. 34 (1981), no. 6, 719--739.


\bibitem[Ve]{Ve} A. Verjovsky, "Arithmetic geometry and dynamics in the unit tangent bundle
of the modular orbifold", in: Dynamical Systems (Santiago 1990), Pitman
Res. Notes Math. No. 285 , Longman Sci. Tech., Harlow, 1993, pp.
263?298


\bibitem[Za1]{Za1}
D. Zagier, Eisenstein Series and the Riemann zeta function, Automorphic Forms, Representation Theory and Arithmetic, Studies in Math. Vol. 10, T.I.F.R., Bombay, 1981, pp. 275-301.


\bibitem[Za2]{Za2} D. Zagier,  The Rankin-Selberg method for authomorphic functions which are not of rapid decay
in J. Fac. Sci. Tokyo 1981.

\bibitem[Za3]{Za3} D. Zagier, Private communication.







\end{thebibliography}
\end{document}